\documentclass{article}
\usepackage{graphicx} 
\usepackage{appendix}
\usepackage{subcaption}
\usepackage{xcolor}
\usepackage{amsmath,amsfonts,latexsym,enumerate}
\usepackage[hyphens]{url}
\usepackage{epsfig,textcomp}
\usepackage{amsthm,longtable,framed,amsbsy,amssymb,amscd}
\usepackage{multirow}
\usepackage{threeparttable}
\usepackage[T1]{fontenc}
\usepackage{palatino,cite}

\usepackage{amssymb,amscd}
\usepackage{nomencl}
\makenomenclature
\makeatletter
\def\myVCENTER#1{\vcenter{\hbox{$\m@th#1$}}}
\makeatother

\long\def\symbolfootnote[#1]#2{\begingroup\def\thefootnote{\fnsymbol{footnote}}\footnote[#1]{#2}\endgroup}

\definecolor{shadecolor}{gray}{0.99}
{\endMakeFramed}

\long\def\symbolfootnote[#1]#2{\begingroup\def\thefootnote{\fnsymbol{footnote}}\footnote[#1]{#2}\endgroup}
\def\qed{\hfill{$\vcenter{\hrule height1pt \hbox{\vrule width1pt height5pt
    \kern5pt \vrule width1pt} \hrule height1pt}$} \medskip}
\newcommand{\m}[1]{{\bf{#1}}}
\newcommand{\g}[1]{\boldsymbol #1}

\newcommand{\C}[1]{{\cal {#1}}}

\title{\bf Margin-Inducing Optimal Guidance \\ and Control for Atmospheric Entry}

\author{Emily M. Palmer\footnote{PhD Candidate, Department of Mechanical and Aerospace Engineering.  E-mail: emilypalmer@ufl.edu.}~ and Anil V. Rao\footnote{Professor, Department of Mechanical and Aerospace Engineering.  Fellow, American Astronautical Society.   Associate Fellow AIAA. \newline\hspace*{12pt} E-mail: anilvrao@ufl.edu.} \vspace{12pt} \\ {\em Department of Mechanical and Aerospace Engineering} \\ {\em University of Florida} \\ {\em Gainesville, FL 32611-6250}}

\date{}

\begin{document}

\maketitle

\nomenclature{$r$}{radial distance from center of planet}
\nomenclature{$v$}{speed, m/s or km/s}
\nomenclature{$\theta$}{longitude, deg or rad}
\nomenclature{$\phi$}{latitude, deg or rad}
\nomenclature{$\gamma$}{flight path angle, deg or rad}
\nomenclature{$\psi$}{Heading angle, rad or deg}
\nomenclature{$\m{c}$}{path constraint}
\nomenclature{$q$}{dynamic pressure, kPa}
\nomenclature{$\dot{Q}$}{heating rate, MW/m$^2$}
\nomenclature{$A$}{sensed acceleration, g}
\nomenclature{$\rho$}{atmospheric density, kg/m$^3$}

\begin{abstract}
A method is presented for margin-inducing optimal guidance and control with application to atmospheric entry.  The method induces margin in state inequality path constraints during entry while simultaneously retaining a large portion of the non-margin-inducing performance.  The margin in the guidance trajectory is based on a differentiable approximation of the maximum of a continuous function. This differentiable approximation is derived from the Kreisselmeier-Steinhauser (KS) function, where the KS function is defined for a function over sampled data. The KS function is first modified into a weighted form where the data points and weights are associated with a multiple-interval Legendre-Gauss-Radau quadrature.  This weighted form is then shown to be an approximation of a logarithmic integral.  This smooth approximation of the maximum function is then used to formulate a margin-inducing optimal control problem, and this margin-inducing optimal control problem is re-solved at specified guidance update times.  The method is applied to a reusable launch vehicle atmospheric entry problem subject to heating rate and dynamic pressure constraints. It is found that the method of this paper produces a guidance trajectory that remains a significant distance from the path constraint limits while maintaining a large portion of the performance of the constrained optimal control problem. 
\end{abstract}

\renewcommand{\baselinestretch}{2}
\normalsize\normalfont

\section{Introduction}\label{section:introduction}

The ability to guide and control a vehicle in the presence of state inequality path constraints is important in many applications.  Such state inequality path constraints include obstacles, keep-out zones, or physical limits that cannot be exceeded.  In applications of this kind, it is not simply undesirable to violate the constraints; constraint violations could lead to disastrous results including crashes, trespassing, or destruction of the moving object. Atmospheric entry is one important application where violating state inequality path constraints could be catastrophic.  Such constraints include heating rate, dynamic pressure, and acceleration load.  In all of these cases, a violation could exceed either the heating capacity or structural load capacity of the vehicle.

Early methods for entry guidance include those used for Apollo (where the vehicle is essentially a drag capsule) and the Shuttle \cite{Harpold1979} (where the vehicle is a lifting body). In the case of Apollo guidance, performance was not the primary objective.  Instead, it was only desired to remain within the thermal and structural capacity of the vehicle while guiding it to a point where a parachute could be deployed.  In the case of Shuttle entry guidance, the goal was again to remain within thermal and structural limits. But in this case, some performance was desired so that the vehicle could maneuver to the terminal area energy management phase and land on a runway.

Because of the limited computational resources on board, these early guidance algorithms were low in performance, which in turn led to low performance guidance trajectories. Specifically, these guidance algorithms employed closed-form guidance laws that required only the execution of simple algebraic operations, allowing for guidance commands to be generated using low-speed on-board processors. In addition, any state inequality path constraints on the vehicle were satisfied indirectly through the conservative control commands that were issued by the guidance algorithm.  

The modern era has seen significant improvement in computational resources, allowing for much higher-speed on-board processors to be utilized for atmospheric flight. Improved flight computers allow for more advanced numerical methods to be executed on-board, increasing the performance that can be obtained. As a result, an entire field known as {\em{computational guidance and control}} (CG\&C)\cite{Lu2017a,tsiotras2017toward} has arisen in recent years, which employs numerical methods to improve performance. In particular, optimization-based CG\&C methods have gained popularity because using an optimization-based CG\&C method often makes it possible to find a feasible solution where otherwise one might not be obtained while simultaneously improving performance.  
 
The development of an optimization-based CG\&C method is centered around solving some type of optimal control or trajectory optimization problem.  A key challenge is that many optimal control problems are nonlinear and do not have analytic solutions.  Moreover, solving such optimal control problems numerically is challenging because it is difficult to find a feasible solution. Even if a feasible solution can in principle be obtained, it is often not possible to do so quickly enough for use on board a vehicle.  Consequently, many optimization-based guidance and control approaches reduce the problem complexity in some way to make the problem more tractable. Key examples of this include model predictive control (MPC) approaches \cite{qin2003survey,mayne2014,eren2017model}, and convex optimization-based methods \cite{wang2024survey,liu2017survey,foust2020optimal,wang2018autonomous}. MPC methods repeatedly re-solve an optimal control problem with the linearized dynamics over a restricted prediction horizon. Recently, nonlinear model predictive control methods have been developed in which the dynamics are not linearized \cite{allgower2000nonlinear}, but the application to complex constrained problems is limited.   In some cases the infeasibility and computational inefficiency can be overcome by the process of {\em convexification}, that is, reformulating the original  non-convex optimal control problem into a convex optimal control problem \cite{mao2016successive,mao2017successive,liu2014solving,szmuk2018successive,blackmore2010minimum}. For some cases, it is possible to apply {\em lossless convexification} \cite{accikmecse2011lossless,harris2014lossless,accikmecse2013lossless} where the convexified problem has the same solution as the original non-convex problem. In other applications, lossless convexification may not be possible.  In such situations, the convexification with loss can be used to obtain solutions with lower performance more reliably than methods for solving the original problem. Additionally, if the convexification of state path constraints is not exact, the solution may not be feasible with respect to the original constraints.

Recently, a computational optimal guidance and control (COG) method was developed in Ref.~\cite{Dennis2019}.  The method of Ref.~\cite{Dennis2019} generates guidance commands by re-solving a nonlinear optimal control problem at the start of each guidance cycle over the remaining horizon. The method transcribes a nonlinear optimal control problem into a general nonlinear programming problem (NLP) using Legendre-Gauss-Radau (LGR) collocation. The method of Ref.~\cite{Dennis2019} has the key feature that it employs a mesh truncation and remapping procedure where the solution from the previous horizon is truncated to leave only the mesh points on the unexpired horizon.  The initial guess is then re-mapped to the mesh on the unexpired horizon. This approach provides the NLP solver with an accurate initial guess when solving the NLP at the start of a guidance cycle. In addition, because the mesh progressively shrinks as the guidance progresses, the NLP also shrinks in size which increases computational efficiency.  

Now, while the method of Ref.~\cite{Dennis2019} does not explicitly prohibit the inclusion of state inequality path constraints, it was not intended for problems where the solutions lie along the path constraint limits.  In particular, the method of Ref.~\cite{Dennis2019} will likely not succeed for any problem where the solution lies along a path constraint limit because uncertainties or disturbances can perturb the solution into a region where the path constraints become infeasible.  Reference~\cite{PalmerRao2022} is an example of a study where the method of Ref.~\cite{Dennis2019} was applied to an atmospheric entry problem with little success because the solution violated state inequality path constraints when simulated with a perturbed model. 

This paper presents a new method for optimal guidance and control that increases margin with respect to state inequality path constraints while maintaining high performance.  Specifically, as opposed to solving a state inequality path-constrained optimal control problem, these path constraints are eliminated from the problem via a margin-inducing objective functional.  As its name implies, the margin-inducing objective functional is designed to increase margin between the state path constraints and their limits. First, a non-differentiable objective functional is developed that is designed to reduce the maximum value of a path constraint function.  Because this objective functional is not differentiable, the terms that involve the maximum function are replaced with smooth approximations. The smooth approximation of the maximum function is derived from the log-sum-exponential Kreisselmeier-Steinhauser (KS) function \cite{kreisselMeier1979} which has been used in constraint aggregation methods\cite{bloss1999dynamic} and penalty methods\cite{wrenn1989indirect,richards2015fast}.  Assuming a monotonically increasing independent variable, a weighted version of the discrete KS function is then derived, where the points are those of a Legendre-Gauss-Radau quadrature approximation of an integral.  The result is that the weighted log-sum-exponential becomes an approximation of a log-integral-exponential.  The terms in the margin-inducing objective functional are then replaced with these differentiable log-integral-exponential terms.  The resulting modified differentiable optimal control problem is then called a {\em margin-inducing optimal control problem (MIOP)}.  A margin-inducing guidance and control method (MIOG) is then developed that includes features from Ref.~\cite{Dennis2019}.  

At first glance, the method of this paper may appear similar to a penalty method. Penalty methods refer to a technique in which a constrained optimal control problem is converted into an unconstrained problem, where the objective functional is augmented with penalty terms on the path constraints \cite{li2011exact,dolgopolik2020exact}. It is noted, however, that for the following reasons, the approach developed in this paper is fundamentally different from a penalty method. In a penalty method, the penalty terms are intended to penalize {\em violations} in the constraints. Typically, these penalty terms are designed such that they are zero when the solution lies within the constraint limits.   Furthermore, the weights on the penalties are chosen such that the solution to the penalized optimal control problem lies in close proximity to the solution of the constrained optimal control problem. Consequently, the solution obtained using a penalty method attempts to recover the optimal value of the original objective functional.  In cases where the optimal solution to the constrained optimal control problem lies along a constraint arc, the solution obtained using a penalty method will also lie along a constraint arc. 

In contrast to a penalty method, the method of this paper is designed to induce margin in the path constraints by reducing the maximum value of the path constraint functions. The resulting solution of the MIOP is designed to lie away from the constraint boundaries. This distinction is key because the modified objective functional developed in this method includes the path constraint functions as a performance metric rather than a penalty. As such, the method of this paper changes the definition of optimality as compared to a penalty method. The result is to produce a guidance trajectory that preserves a large portion of the original performance while providing the key benefit that margin with respect to the path constraints is maintained in the presence of uncertainties or environmental disturbances. 

It is noted that the method of this paper builds upon the work of Ref.~\cite{PalmerRaoCCG2026}, in which the original inequality path constraints are explicitly enforced, and objective weights are switched based on the proximity of the guidance trajectory to the constraint limits.  On the other hand, in this paper the state inequality path constraints are removed, and margin-inducing terms are included in the objective functional and a single optimal control problem formulation is used throughout guidance. Consequently, margin is promoted through the objective functional without explicitly enforcing the state path constraints. Finally, this paper provides a mathematical justification that the log-integral function is an approximation for the maximum of a continuous function.

The contributions of this research are as follows.  First, a path-constrained optimal control problem is reformulated by eliminating the path constraints and designing a modified objective functional that is intended to increase the margin between the path constraint functions and their limits. This objective functional, however, is not differentiable because it employs a maximum function.  In order to make the objective functional differentiable, a smooth approximation of the maximum function is derived using a continuous analogue of the Kreisselmeier-Steinhauser function \cite{kreisselMeier1979}.  This differentiable approximation of the maximum function is then used to replace the non-differentiable terms in the aforementioned modified objective functional, creating a {\em margin-inducing objective functional} and an associated {\em margin-inducing optimal control problem} (MIOP). This MIOP is approximated using Legendre-Gauss-Radau (LGR) collocation, and a margin-inducing guidance and control (MIOG) method is developed by re-solving this MIOP using features from the method of Ref.~\cite{Dennis2019}.  The MIOG method has the key feature that it produces a guidance trajectory with desired qualities that retains a large portion of the performance that would have been obtained by optimizing the original constrained optimal control problem. The MIOG method is then demonstrated on a reusable launch vehicle (RLV) entry problem subject to heating rate and dynamic pressure path constraints.  It is found that the MIOG method is capable of re-solving the MIOP at the guidance update times quickly while controlling the vehicle so that it remains feasible with respect to the path constraints.   

This paper is organized as follows.  Section~\ref{sect:background} provides background information.  Section~\ref{sect:miop} provides the derivation of the margin-inducing optimal control problem that forms the basis of the developed guidance and control method.  Section~\ref{sect:miog} describes the margin-inducing guidance and control method.  Section~\ref{sect:rlv-ocp} describes the state inequality path-constrained reusable launch vehicle entry problem to which the method developed in Section~\ref{sect:miog} is applied.  Section~\ref{sect:rlv-results-miog} provides the results of applying the developed method to the reusable launch vehicle entry problem. Section~\ref{sect:summary} provides a summary of the key results obtained in this study.  Finally, Section~\ref{sect:conclusions} provides conclusions.  

 \section{Background Information}\label{sect:background}

This section provides the important background information required to develop the margin-inducing optimal guidance (MIOG) method.  First, two dynamic models of interest are presented: (1) a reference dynamic model and (2) an actual dynamic model. The reference dynamic model is the one used to solve an optimal control problem, while the actual dynamic model represents the dynamics corresponding to the actual motion of the system.  Second, the state path-constrained Bolza optimal control problem relevant to this paper is formulated. Third, the state path-constrained optimal control problem is discretized using multiple-interval Legendre-Gauss-Radau (LGR) collocation.  This background information is used in the remainder of the paper both to develop the MIOG method and to demonstrate the MIOG method described in Section \ref{sect:miog} on the path-constrained atmospheric entry problem described in Section \ref{sect:rlv-ocp}. 

\subsection{Reference and Actual Dynamics}\label{sect:dynamic-models}

Two dynamic models are of interest,  the {\em reference} dynamics and the {\em actual} dynamics.  The reference dynamics are those upon which the method for guidance and control developed in this paper is based. The {\em actual} dynamics of the system are those that are realized by the system during the actual motion. Each of these models is described below.

\subsubsection{Reference Dynamics}\label{sect:reference-dynamic-model}

The reference dynamics are divided into {\em state dynamics} 
\begin{equation}\label{eq:reference-state-dynamics}
\begin{array}{lcl}
  \displaystyle \frac{d\m{x}}{dt}=\m{f}(\m{x}(t),\m{u}(t)), \\
\end{array}
\end{equation}
and control (actuator) dynamics
\begin{equation}\label{eq:reference-control-dynamics}
\begin{array}{lcl}
  \displaystyle \frac{d\m{u}}{dt}=\m{a}(\m{u}(t),\m{p}(t)), \\
\end{array}
\end{equation}
where $\m{x}(t)\in\mathbb{R}^{n_x}$ is the state, $\m{u}(t)\in\mathbb{R}^{n_u}$ is the control, $\m{p}(t)\in\mathbb{R}^{n_u}$ is the pseudo-control, $\m{f}: \mathbb{R}^{n_x} \times \mathbb{R}^{n_u} \rightarrow \mathbb{R}^{n_x}$ is the vector field for the reference state dynamics, $\m{a}: \mathbb{R}^{n_u} \times \mathbb{R}^{n_u} \rightarrow \mathbb{R}^{n_u}$ is the vector field for the reference control (actuator) dynamics. In this study, all optimal control problems will be formulated in terms of the reference dynamics $(\m{f},\m{a})$.  Furthermore, with regard to the optimal control problem the "state" will be $(\m{x}(t),\m{u}(t))$ because the actuator dynamics are part of the reference model.  

\subsubsection{Actual Dynamics}\label{sect:actual-dynamic-model}

The actual state dynamics, that is, the dynamics that correspond to the actual realized motion, are given as
\begin{equation}\label{eq:actual-dynamics}
  \frac{d\m{\tilde{x}}}{dt}=\m{\tilde{f}}(\m{\tilde{x}}(t),\m{u}(t)), 
\end{equation}
where $\tilde{\m{f}}:\mathbb{R}^{n_x} \times \mathbb{R}^{n_u}  \rightarrow \mathbb{R}^{n_x}$, and $\m{\tilde{x}}(t) \in \mathbb{R}^{n_x}$ is the state of the actual system. In this study, the actual dynamic model defined by the vector field $\tilde{\m{f}}$ is used to simulate the motion of the actual system.  Then, because the control, $\m{u}(t)$, is the input to the actual state dynamics in Eq.~\eqref{eq:actual-dynamics}, the control (actuator) dynamics are omitted. Due to disturbances, measurement errors, or other uncertainties, the reference dynamics will differ from the actual dynamics (that is, $\m{f} \neq \m{\tilde{f}})$. 

\subsection{Constrained Optimal Control Problem}\label{sect:constrained-ocp}

Consider the following state path-constrained optimal control problem defined on $t\in[t_0,t_f]$, where $t_0$ and $t_f$ are the initial and terminal times, respectively. Determine the state, $\m{x}(t)\in\mathbb{R}^{n_x}$, the control, $\m{u}(t)\in\mathbb{R}^{n_u}$, the pseudo-control, $\m{p}(t)\in\mathbb{R}^{n_u}$, and the terminal time $t_f$ that minimize the objective functional 
\begin{equation}\label{eq:obj}
\C{J} = \C{M}(\m{x}(t_f),t_f) + \int_{t_0}^{t_f} \C{L}(\m{x}(t),\m{u}(t),\m{p}(t))\, dt,
\end{equation}
subject to the state dynamic constraints
\begin{equation}\label{eq:dyncon}
  \frac{d\m{x}}{dt}-\m{f}(\m{x}(t),\m{u}(t))=\m{0},
\end{equation}
the control (actuator) dynamic constraints, 
\begin{equation}\label{eq:actcon}
  \frac{d\m{u}}{dt}-\m{a}(\m{u}(t),\m{p}(t))=\m{0},
\end{equation}
the control constraints
\begin{equation}\label{eq:control-constraints}
\begin{array}{lcccl}
    \m{u}_{\min} &\leq& \m{u}(t) &\leq &\m{u}_{\max},
\end{array}
\end{equation}
the pseudo-control constraints
\begin{equation}\label{eq:pseudo-control-constraints}
\begin{array}{lcccl}
    \m{p}_{\min} &\leq& \m{p}(t) &\leq &\m{p}_{\max},
\end{array}
\end{equation}
the boundary conditions 
\begin{equation}\label{eq:boundcond}
\begin{array}{rcl}
\m{x}(t_0) &=& \m{x}_0 = \textrm{Fixed}, \\
 \m{b}(\m{x}(t_f),t_f) &=& \m{0}, \\
t_0 &= &\mathrm{Fixed},
\end{array}
\end{equation}
and the state path constraints
\begin{equation}\label{eq:pathcon}
    \m{c}(\m{x}(t)) \leq \m{c}_{\max}. 
\end{equation}
The function $\C{M}: \mathbb{R}^{n_x} \times \mathbb{R} \rightarrow \mathbb{R}$ is the Mayer cost, $\C{L}: \mathbb{R}^{n_x} \times \mathbb{R}^{n_u} \times \mathbb{R}^{n_u} \rightarrow \mathbb{R}$ is the Lagrangian, $\m{b}: \mathbb{R}^{n_x} \times \mathbb{R} \rightarrow \mathbb{R}^{n_b}$ is the function that describes the $n_b$ terminal constraints, $\m{c}:\mathbb{R}^{n_x}\rightarrow \mathbb{R}^{n_c}$ describes the $n_c$ state path constraints, and $\m{x}_0\in\mathbb{R}^{n_x}$ is the initial state.  

Finally, because it will be useful in formulating the LGR collocation method described in Section \ref{sect:lgr-nlp}, the optimal control problem can be transformed from the time interval $t \in [t_0,t_f]$ to the interval $\tau\in[-1,+1]$ via the transformation
\begin{equation}\label{eq:tau-to-t}
t = \frac{t_f-t_0}{2}\tau + \frac{t_f+t_0}{2}. 
\end{equation}
The constrained optimal control problem is then given in terms of the variable $\tau\in[-1,+1]$ as follows.  Minimize the objective functional
\begin{equation}\label{eq:obj-tau}
    \C{J} = \C{M}(\m{x}(+1),t_f) + \frac{t_f-t_0}{2}\int_{-1}^{+1} \C{L}(\m{x}(\tau),\m{u}(\tau),\m{p}(\tau))\, d\tau
\end{equation}
subject to the state dynamic constraints
\begin{equation}\label{eq:dyncon-tau}
  \frac{d\m{x}}{d\tau}-\frac{t_f-t_0}{2}\m{f}(\m{x}(\tau),\m{u}(\tau))=\m{0}, 
\end{equation}
the control (actuator) dynamic constraints
\begin{equation}\label{eq:actcon-tau}
  \frac{d\m{u}}{d\tau}-\frac{t_f-t_0}{2}\m{a}(\m{u}(\tau),\m{p}(\tau))=\m{0}, 
\end{equation}
the boundary conditions
\begin{equation}\label{eq:boundcond-tau}
\begin{array}{rcl}
\m{x}(-1) &=& \m{x}_0 = \textrm{Fixed}, \\
 \m{b}(\m{x}(+1),t_f) &=& \m{0}, \\
  t_0 &= &\mathrm{Fixed},
\end{array}
\end{equation}
the control constraints
\begin{equation}\label{eq:control-constraints-tau}
    \m{u}_{\min} \leq \m{u}(\tau) \leq \m{u}_{\max},
\end{equation}
the pseudo-control constraints
\begin{equation}\label{eq:pseudo-control-constraints-tau}
    \m{p}_{\min} \leq \m{p}(\tau) \leq \m{p}_{\max},
\end{equation}
and the state path constraints
\begin{equation}\label{eq:pathcon-tau}
    \m{c}(\m{x}(\tau)) \leq \m{c}_{\max}. 
\end{equation}

\subsection{Legendre-Gauss-Radau Collocation for Optimal Control}\label{sect:lgr-nlp}

In this paper all optimal control problems are solved numerically using the multiple-interval $hp$ Legendre-Gauss-Radau (LGR) collocation method as described in Refs.~\cite{Garg2010,Garg2011b,Hager2015a}, where multiple-interval LGR collocation transcribes the optimal control problem into a large sparse nonlinear programming problem (NLP). In this section the constrained optimal control problem on $\tau\in[-1,+1]$ described in Section \ref{sect:constrained-ocp} is approximated using multiple-interval LGR collocation.   

Suppose that the interval $\tau\in[-1,+1]$ is partitioned into $K$ subintervals $\C{I}_k\;(k=1,\ldots,K)$, where $\C{I}_k=[T_{k-1},T_k]$, and $-1=T_0<T_1<\cdots<T_K=+1$.  Let $\tau^{(k)}=(\tau_1^{(k)},\ldots,\tau_{N_k}^{(k)})$ and $\m{w}^{(k)}=(w_1^{(k)},\ldots,w_{N_k}^{(k)})$ be the $N_k$ Legendre-Gauss-Radau points and weights, respectively, defined on mesh interval $\C{I}_k,\; k\in[1,\ldots,K]$.  Furthermore, let $\g{\tau}=(\g{\tau}^{(1)},\ldots,\g{\tau}^{(K)})$ and $\m{w}=(\m{w}^{(1)},\ldots,\m{w}^{(K)})$ be the cumulative LGR points and weights on $\tau\in[-1,+1]$, where $N=\sum_{k=1}^{K} N_k$. It is noted that $\g{\tau}$ is a vector with monotonically increasing values, that is,  $\tau_1<\tau_2<\cdots<\tau_N$.  Finally, let $\tau_{N+1}=+1$.  Then, if $(\m{X}_i,\m{U}_i,\m{P}_i)$ are the approximations of the state, control, and pseudo-control at $\tau=\tau_i$ and $(\m{X}_{N+1},\m{U}_{N+1})$ are the approximations of the state and control at $\tau_{N+1}=+1$, the constrained optimal control problem is approximated using multiple-interval LGR collocation as follows.  Minimize the objective function
\begin{equation}\label{eq:lgr-obj}
\C{J} \approx J  = \C{M}(\m{X}_{N+1},t_f) + \frac{t_f-t_0}{2}\sum_{i=1}^N w_i\C{L}\left[\m{X}_i,\m{U}_i,\m{P}_i\right]
\end{equation}
subject to the collocated state dynamic constraints
\begin{equation}\label{eq:lgr-state-dynamics}
    \sum_{j=1}^{N+1}D_{ij}\m{X}_j-\frac{t_f-t_0}{2}\m{f}(\m{X}_i,\m{U}_i) = \m{0}, \quad (i=1,\ldots,N),
\end{equation}
the collocated control (actuator) dynamic constraints
\begin{equation}\label{eq:lgr-control-dynamics}
    \sum_{j=1}^{N+1}D_{ij}\m{U}_j-\frac{t_f-t_0}{2}\m{a}(\m{U}_i,\m{P}_i) = \m{0}, \quad (i=1,\ldots,N),
\end{equation}
the discretized boundary conditions
\begin{equation}\label{eq:lgr-bc}
\begin{array}{rcl}
\m{X}_1 &=& \m{x}_0, \\
\m{b}(\m{X}_{N+1},t_f) &=& \m{0}, \\
t_0 &= &\mathrm{Fixed},
\end{array}
\end{equation}
the discretized control constraints
\begin{equation}\label{eq:lgr-controlcon}
    \m{u}_{\min} \leq \m{U}_i \leq \m{u}_{\max},\quad (i=1,\ldots,N),
\end{equation}
the discretized pseudo-control constraints
\begin{equation}\label{eq:lgr-pseudo-controlcon}
    \m{p}_{\min} \leq \m{P}_i \leq \m{p}_{\max},\quad (i=1,\ldots,N),
\end{equation}
and the discretized state path constraints
\begin{equation}\label{eq:lgr-pathcon}
    \m{c}(\m{X}_i) \leq \m{c}_{\max},\quad (i=1,\ldots,N).
\end{equation}
Finally, it is noted in Eq.~\eqref{eq:lgr-state-dynamics} that $D_{ij}$ are the elements of the $N\times(N+1)$ LGR differentiation matrix on $\tau\in[-1,+1]$, where the LGR differentiation matrix has a block-diagonal structure \cite{PattersonRao2015b}, such that the $k^{th}$ block is used to compute derivative approximations on mesh interval $\C{I}_k$.  Further details on the structure of the LGR NLP can be found in Ref.~\cite{PattersonRao2015b}.  

\section{Margin-Inducing Optimal Control Problem}\label{sect:miop}

This section describes the margin-inducing optimal control problem (MIOP) that is used as the basis for the margin-inducing optimal guidance (MIOG) method of Section~\ref{sect:miog}.  The MIOP described in this section removes the state path constraints and replaces them with terms in the objective functional intended to improve the margin between the path constraints and their limits. Because these additional terms involve the non-differentiable maximum function, they are replaced with a smooth approximation of the maximum function.  This smooth approximation is developed by starting with the Kreisselmeier-Steinhauser (KS) function, where the KS function is used in this paper to approximate the maximum of a function that depends upon a monotonically increasing discrete independent variable.  Specifically, the KS function is replaced with a weighted modification of the KS function, where the weights are those of a multiple-interval LGR quadrature as described in Section \ref{sect:lgr-nlp}. This multiple-interval LGR quadrature approximation is then used as an approximation of the logarithm of an integral (log-integral), and the terms involving the maximum function in the margin-inducing objective functional are replaced with log-integral terms. Finally, the MIOP is stated formally using this differentiable form of the margin-inducing objective functional.  

\subsection{Use of Maximum Function to Increase Path Constraint Margin}

The definition of margin in a path constraint is the difference between the path constraint limit and the maximum value of the path constraint function (where positive margin means the constraint is satisfied while negative margin means the constraint is violated).  In other words, given a single state path constraint of the form given by Eq.~\eqref{eq:pathcon} the margin is defined as
\begin{equation}\label{eq:definition-of-margin}
  c_{\max} - \max_{t\in[t_0,t_f]} c(\m{x}(t)) = c_{\max}\left[1-\max_{t\in[t_0,t_f]}\frac{c(\m{x}(t))}{c_{\max}}\right],
\end{equation}
where it is assumed that $c_{\max}>0$.  Then, increasing margin is equivalent to decreasing
\begin{equation}\label{eq:reduce-max}
    \max_{t\in[t_0,t_f]} \frac{c(\m{x}(t))}{c_{\max}}.
\end{equation}
Suppose that with the goal of improving path constraint margin, the state path constraints given in Eq.~\eqref{eq:pathcon} are eliminated from the optimal control problem and terms similar to that shown in Eq.~\eqref{eq:reduce-max} are included in the objective functional. Then, additional terms in the objective functional that increase path constraint margin would have the form
\begin{equation}\label{eq:augumented-terms-objective-max}
    \C{J}_I = \sum_{j=1}^{n_c} \max_{t\in[t_0,t_f]}\left(e_j\right),
\end{equation}
where
\begin{equation}
    e_j = \frac{c_j(\m{x}(t))}{c_{\max,j}}, \quad (c_{\max,j}>0).
\end{equation}
Then a margin-inducing objective functional would have the form 
\begin{equation}\label{eq:augmented-objective-max}
    \C{J}_a = \C{J}_{\textrm{ref}} + \C{J}_I = \C{M} + \int_{t_0}^{t_f} \C{L}~dt + \sum_{j=1}^{n_c}  \max_{t\in[t_0,t_f]}\left(e_j\right).
\end{equation}
While in principle the objective functional given in Eq.~\eqref{eq:augmented-objective-max} could be used to induce margin in the solution of the optimal control problem, the non-differentiability of the maximum function leads to a non-differentiable optimal control problem formulation which is unsuitable when using a gradient-based method such as those employed in this paper.  In Section~\ref{sect:smooth-max-approximation} a differentiable approximation of the maximum function is derived which is then used to modify Eq.~\eqref{eq:augmented-objective-max}. 

\subsection{Differentiable Approximation of Maximum Function}\label{sect:smooth-max-approximation}

In this paper the non-differentiable objective functional in Eq.~\eqref{eq:augmented-objective-max} is approximated with a differentiable objective functional in terms of a logarithm of an integral.  Consequently, it is necessary to begin the derivation with the assumption that the function whose maximum is desired depends upon a monotonically increasing variable. To this end, let $f(\tau)$ be a differentiable function on $\tau\in[-1,+1]$, where $f(\tau)$ is sampled at the monotonically increasing $N$ LGR points $\tau_i,\;(i=1,\ldots,N)$ as defined in Section \ref{sect:lgr-nlp}.  Then the maximum of $f(\tau)$ over these $N$ LGR points is given as
\begin{equation}\label{eq:max-sampled-data}
    f_{\max} = \max_{i\in[1,\ldots,N]} f(\tau_i).
\end{equation}
It is seen that the function in Eq.~\eqref{eq:max-sampled-data} is not differentiable.  Therefore, Eq.~\eqref{eq:max-sampled-data} cannot be included as part of a smooth formulation of an optimal control problem and a smooth approximation must be developed.  The remainder of this section derives a smooth approximation of Eq.~\eqref{eq:max-sampled-data}.

It has been shown that the non-differentiable maximum function in Eq.~\eqref{eq:max-sampled-data} can be approximated using the following differentiable {\em Kreisselmeier-Steinhauser} (KS) function \cite{kreisselMeier1979}:
\begin{equation}\label{eq:kreisselmeier-steinhauser-sampled-data}
f_{\max} = \max_{i\in[1,\ldots,N]} f(\tau_i) \approx KS(f) = \frac{1}{\beta}\log\left\{\sum_{i=1}^{N} \exp\left[\beta f(\tau_i)\right]\right\},
\end{equation}
where $\log(\cdot)$ is the natural logarithm and the parameter $\beta$ controls the accuracy of the approximation.  Suppose now that the function in Eq.~\eqref{eq:kreisselmeier-steinhauser-sampled-data} is replaced with the weighted form
\begin{equation}\label{eq:kreisselmeier-steinhauser-sampled-data-weighted}
f_{\max} = \max_{i\in[1,\ldots,N]} f(\tau_i) \approx KS(f) = \frac{1}{\beta}\log\left\{\sum_{i=1}^{N} w_i \exp\left[\beta f(\tau_i)\right]\right\},
\end{equation}
where 
$w_i\;(i=1,\ldots,N)$ are the $N$ LGR quadrature weights corresponding to the $N$ LGR points $\tau_i,\;(i=1\ldots,N)$. It is then seen that the sum in Eq.~\eqref{eq:kreisselmeier-steinhauser-sampled-data-weighted} is an LGR approximation of the integral of $\exp\left[\beta f(\tau)\right]$ on $\tau\in[-1,+1]$, that is, 
\begin{equation}
\int_{-1}^{+1} \exp\left[\beta f(\tau)\right]d\tau \approx \sum_{i=1}^{N} w_i \exp\left[\beta f(\tau_i)\right].   
\end{equation}
Consequently, Eq.~\eqref{eq:kreisselmeier-steinhauser-sampled-data-weighted} is an approximation of the following logarithmic integral: 
\begin{equation}\label{eq:kreisselmeier-steinhauser-continuous-from-weighted-discrete-tau}
    \frac{1}{\beta}\log\left[\int_{-1}^{+1} \exp\left[\beta f(\tau)\right] d\tau\right] \approx \frac{1}{\beta}\log\left\{\sum_{i=1}^{N} w_i \exp\left[\beta f(\tau_i)\right]\right\}.
\end{equation}
Adding the term $\log\left[(t_f-t_0)/2\right]/\beta$ to both sides of Eq.~\eqref{eq:kreisselmeier-steinhauser-continuous-from-weighted-discrete-tau} then gives
\begin{equation}\label{eq:kreisselmeier-steinhauser-continuous-from-weighted-discrete-pre}
    \frac{1}{\beta}\log\left[\frac{t_f-t_0}{2}\int_{-1}^{+1} \exp\left[\beta f(\tau)\right] d\tau\right] \approx \frac{1}{\beta}\log\left\{\frac{t_f-t_0}{2} \sum_{i=1}^{N} w_i \exp\left[\beta f(\tau_i)\right]\right\}.
\end{equation}
Comparing the term inside the logarithm on the left-hand side of Eq.~\eqref{eq:kreisselmeier-steinhauser-continuous-from-weighted-discrete-pre} with Eq.~\eqref{eq:obj-tau}, it is seen that
\begin{equation}\label{eq:integral-tau-and-t}
\frac{t_f-t_0}{2}\int_{-1}^{+1} \exp\left[\beta f(\tau)\right] d\tau = \int_{t_0}^{t_f} \exp\left[\beta f(\tau(t))\right] dt \equiv  \int_{t_0}^{t_f} \exp\left[\beta f(t)\right] dt, 
\end{equation}
where, using the inverse of the transformation in Eq.~\eqref{eq:tau-to-t}, $f(t)$ is a shorthand for $f(\tau(t))$.  Then, define the function $\C{E}_{\beta}$ as
\begin{equation}\label{eq:kreisselmeier-steinhauser-continuous-from-weighted-discrete}
    \C{E}_{\beta}(f) = \frac{1}{\beta}\log\left[\int_{t_0}^{t_f} \exp\left[\beta f(t)\right] dt\right] \approx \frac{1}{\beta}\log\left\{\frac{t_f-t_0}{2} \sum_{i=1}^{N} w_i \exp\left[\beta f(\tau_i)\right]\right\}.
\end{equation}
The maximum of a function $f(t)$ on $t\in[t_0,t_f]$ is then approximated as
\begin{equation}\label{eq:max-approximation-continuous}
    \max_{t\in[t_0,t_f]} f(t) \approx \C{E}_{\beta}(f) = 
    \frac{1}{\beta}\log\left[\int_{t_0}^{t_f} \exp\left[\beta f(t)\right] dt\right],
\end{equation}
where $\beta>0$ and Eq.~\eqref{eq:max-approximation-continuous} becomes exact in the limit as $\beta\rightarrow\infty$ via the Laplace principle \cite{varadhan1966asymptotic,dembozeitouni1998}.  The function $\C{E}_{\beta}(\cdot)$ given in Eq.~\eqref{eq:max-approximation-continuous} is used in the remainder of this paper to develop the margin-inducing objective functional described in Section \ref{sect:margin-inducing-objective} and the associated margin-inducing optimal control problem (MIOP) described in Section \ref{sect:margin-inducing-ocp}. 

\begin{figure}[htbp]
    \centering

    \begin{minipage}[t]{0.48\textwidth}
        \begin{subfigure}[t]{\linewidth}
            \centering
            \includegraphics[width=\linewidth]{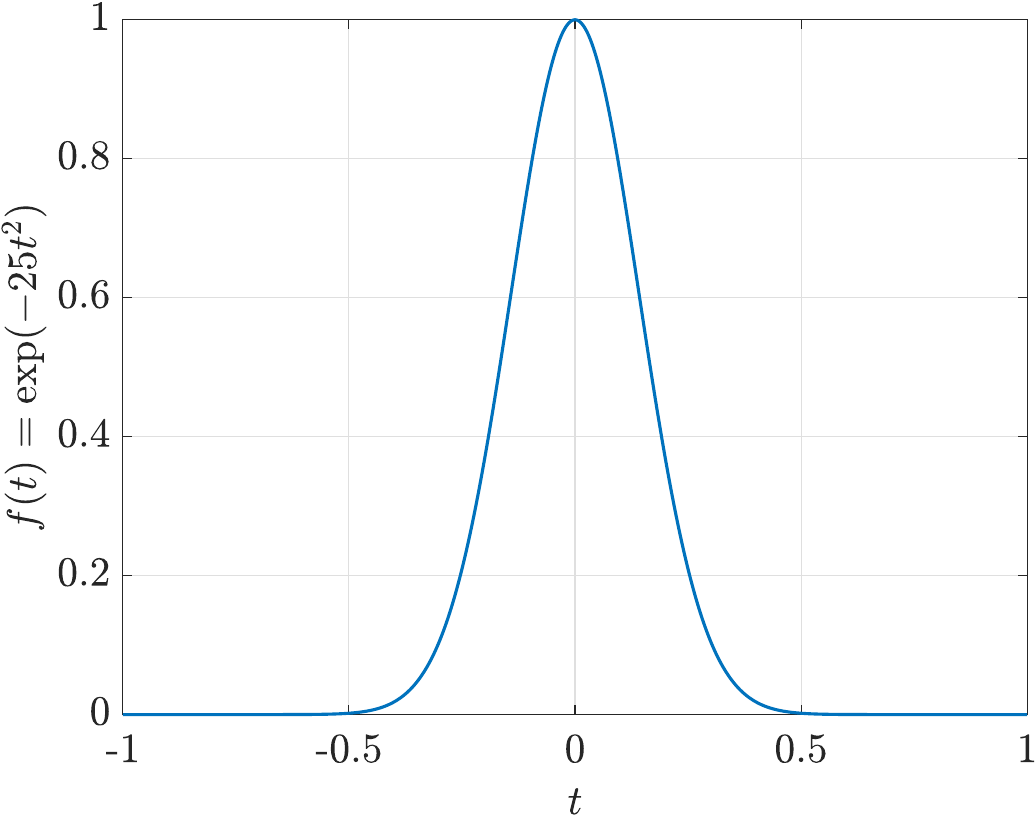}
            \caption{$f(t)=\exp(-25t^2)$.}
            \label{fig:example-Ebeta-a}
        \end{subfigure}
    \end{minipage}
    \hfill
    \begin{minipage}[t]{0.48\textwidth}
        \begin{subfigure}[t]{\linewidth}
            \centering
            \includegraphics[width=\linewidth]{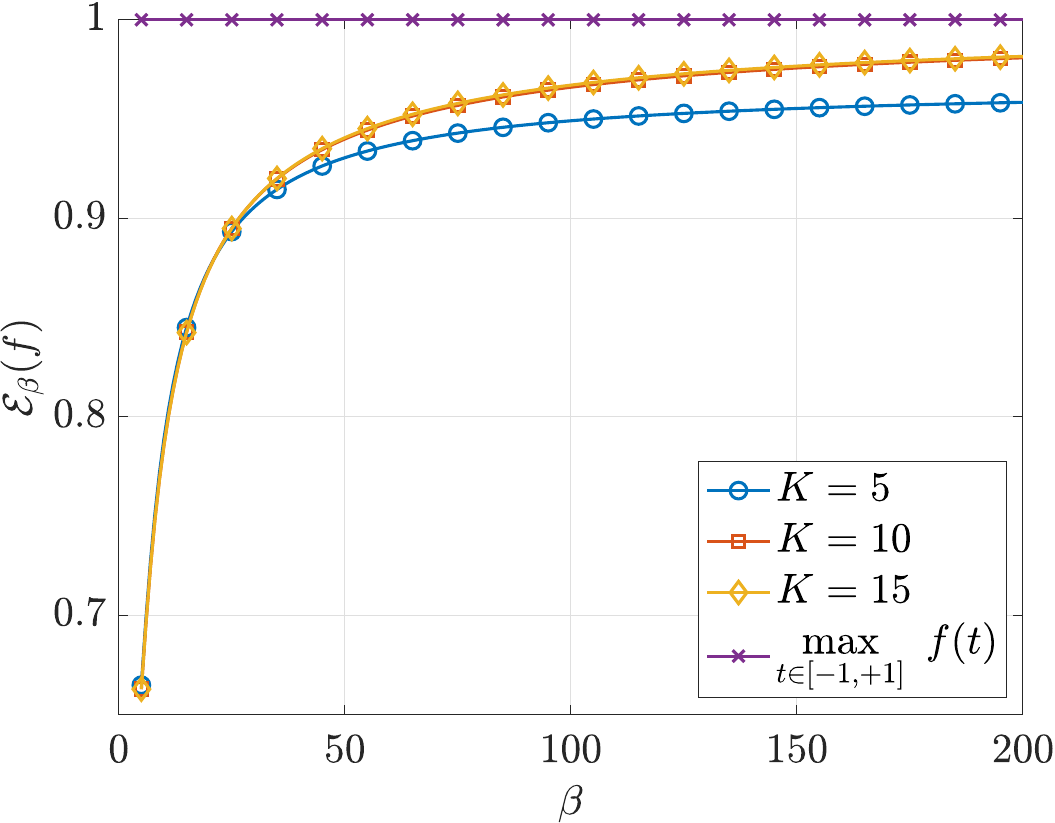}
            \caption{Approximation of smooth maximum, $\C{E}_{\beta}(f)$}
            \label{fig:example-Ebeta-b}
        \end{subfigure}
    \end{minipage}

    \caption{$\C{E}_{\beta}(f)$ for $f(t)=\exp(-25t^2)$ and $t\in[-1,+1]$ using a four-point-per-interval LGR quadrature with $K=[5,10,15]$ equally-spaced intervals.}
    \label{fig:example-Ebeta}
\end{figure}

To see how the approximation in Eq.~\eqref{eq:max-approximation-continuous} works in practice, consider the function $f(t)=\exp(-25t^2)$ on $t\in[-1,+1]$.  Figure~\ref{fig:example-Ebeta-a} shows $f(t)$ while Fig.~\ref{fig:example-Ebeta-b} shows the approximation of the maximum of $f(t)$ given in
Eq.~\eqref{eq:max-approximation-continuous}
for $\beta\in[5,200]$ using an LGR quadrature with four LGR points per interval and $K=[5,10,15]$ equally-spaced intervals. It is seen for this example that the approximation $\C{E}_{\beta}(f)$ of the maximum of $f(t)$ improves as $\beta$ increases.  In practice, while Eq.~\eqref{eq:max-approximation-continuous} may be less accurate for smaller values of $\beta$, smaller values of $\beta$ are preferable because the integral in Eq.~\eqref{eq:max-approximation-continuous} may become unreasonably large as $\beta$ increases.

\subsection{Differentiable Margin-Inducing Objective Functional}\label{sect:margin-inducing-objective}

Suppose now that the terms $\displaystyle\max_{t\in[t_0,t_f]}\left(e_j\right),\;(j=1,\ldots,n_c)$ in the objective functional of Eq.~\eqref{eq:augmented-objective-max} are replaced with approximations based on Eq.~\eqref{eq:max-approximation-continuous} as follows:
\begin{equation}\label{eq:max-terms-in-augmented-objective}
\max_{t\in[t_0,t_f]}\left(e_j\right)\approx \C{E}_{\beta}(e_j) = \frac{1}{\beta}\log\left[\int_{t_0}^{t_f} \exp\left[\beta e_j\right] dt\right],\quad (j=1,\ldots,n_c). 
\end{equation}
Substituting the result of Eq.~\eqref{eq:max-terms-in-augmented-objective} into Eq.~\eqref{eq:augmented-objective-max} then gives
\begin{equation}\label{eq:margin-inducing-objective}
    \C{J}_a = \C{J}_{\textrm{ref}} + \C{J}_I = \C{M} + \int_{t_0}^{t_f} \C{L}~dt + \sum_{j=1}^{n_c} \C{E}_{\beta}(e_j). 
\end{equation}
The objective functional defined in Eq.~\eqref{eq:margin-inducing-objective} is called a {\em margin-inducing objective functional}.  It is seen that the margin-inducing objective functional in Eq.~\eqref{eq:margin-inducing-objective} is designed to reduce the maximum value of each of the path constraint functions, thereby increasing the margin between the path constraint functions and their limits.

\subsection{Margin-Inducing Optimal Control Problem}\label{sect:margin-inducing-ocp}

The margin-inducing optimal control problem (MIOP) is stated as follows.  First, remove all inequality path constraints given in Eq.~\eqref{eq:pathcon}.  Then the margin-inducing optimal control problem is stated as follows.  Determine the state, $\m{x}(t)\in\mathbb{R}^{n_x}$, the control, $\m{u}(t)\in \mathbb{R}^{n_u}$, the pseudo-control, $\m{p}(t) \in \mathbb{R}^{n_u}$, and the terminal time $t_f$ on $t\in[t_0,t_f]$ that minimize the objective functional 
\begin{equation}\label{eq:objective-miop}
    \C{J}_a = \C{M}(\m{x}(t_f),t_f) + \int_{t_0}^{t_f} \C{L}(\m{x}(t),\m{u}(t),\m{p}(t))\, dt +  \sum_{j=1}^{n_c} \C{E}_{\beta}(e_j). 
\end{equation}
subject to the dynamic constraints, control (actuator) dynamic constraints, control constraints, pseudo-control constraints, and boundary conditions given in Eqs.~\eqref{eq:dyncon}--\eqref{eq:boundcond}, respectively.  It is stated again that the MIOP does not contain the state path constraints given in Eq.~\eqref{eq:pathcon} because these state path constraints have been incorporated into the formulation of the MIOP objective functional as given in Eq.~\eqref{eq:objective-miop}.  In this paper the MIOP is solved using the LGR collocation method described in Section \ref{sect:lgr-nlp}.  

\section{Guidance Methodology}\label{sect:miog}

This section develops the margin-inducing optimal guidance (MIOG) method.  The MIOG method re-solves the MIOP described in Section \ref{sect:miop} at the start of each guidance cycle. The control obtained by re-solving the MIOP is then used in the actual dynamics on the current guidance cycle.  Furthermore, if the LGR NLP corresponding to the MIOP is infeasible for any set of consecutive guidance cycles, the control obtained from the last solution obtained to the LGR NLP is used until the next LGR NLP solution is obtained.  The process of re-solving the optimal control problem and simulating the motion is repeated until the remaining time horizon is less than the guidance cycle duration.  

\subsection{Margin-Inducing Optimal Guidance}



Consider re-solving the MIOP at guidance cycles $[t_0^{(s)},t_e^{(s)}]$, where
\begin{eqnarray}\label{eq:guidance-update-start-end}
    t_0^{(s)} & = & t_0 + s \Delta T, \\
    t_e^{(s)} & = & t_0 + (s+1)\Delta T, 
\end{eqnarray}
are the start and terminus times of guidance cycle $s\in[1,\ldots,S]$, $s=0$ corresponds to the initial time $t=t_0$, $\Delta T$ is the guidance cycle duration, and $S$ is the total number of guidance cycles.  The MIOP is then re-solved at the start of each guidance cycle $s\in[1,\ldots,S]$ (that is, every $\Delta T$ time units) on the remaining horizon $t\in[t_0^{(s)},t_f]$ with the initial state $\m{x}(t_0^{(s)})=\tilde{\m{x}}(t_e^{(s-1)})$, where $\tilde{\m{x}}(t_e^{(s-1)})$ is the state of the actual system (that is, the system corresponding to the actual dynamics as modeled in Eq.~\eqref{eq:actual-dynamics}). 

\begin{figure}[h]
    \centering    \includegraphics[width=0.9\linewidth]{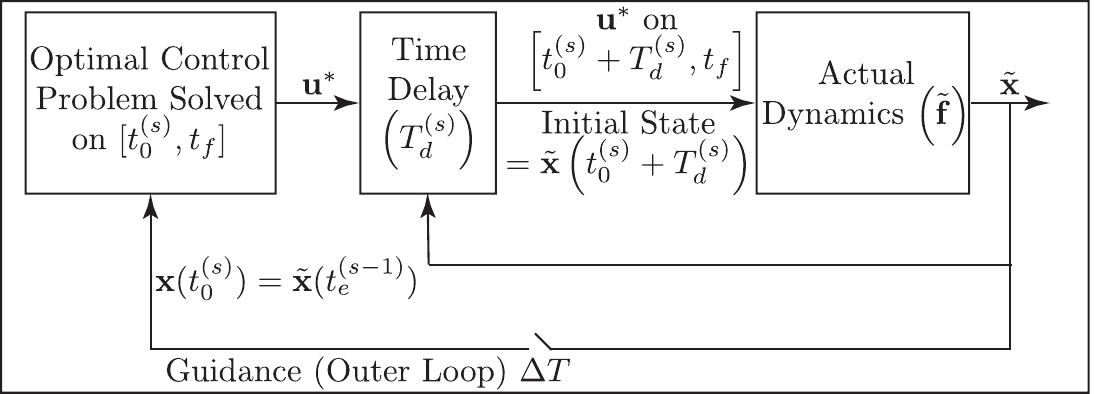}
    \caption{Block diagram of the margin-inducing optimal guidance (MIOG) method.}
    \label{fig:blockdiagram}
\end{figure}

\begin{figure}[h]
\centering
\includegraphics[scale=0.75]{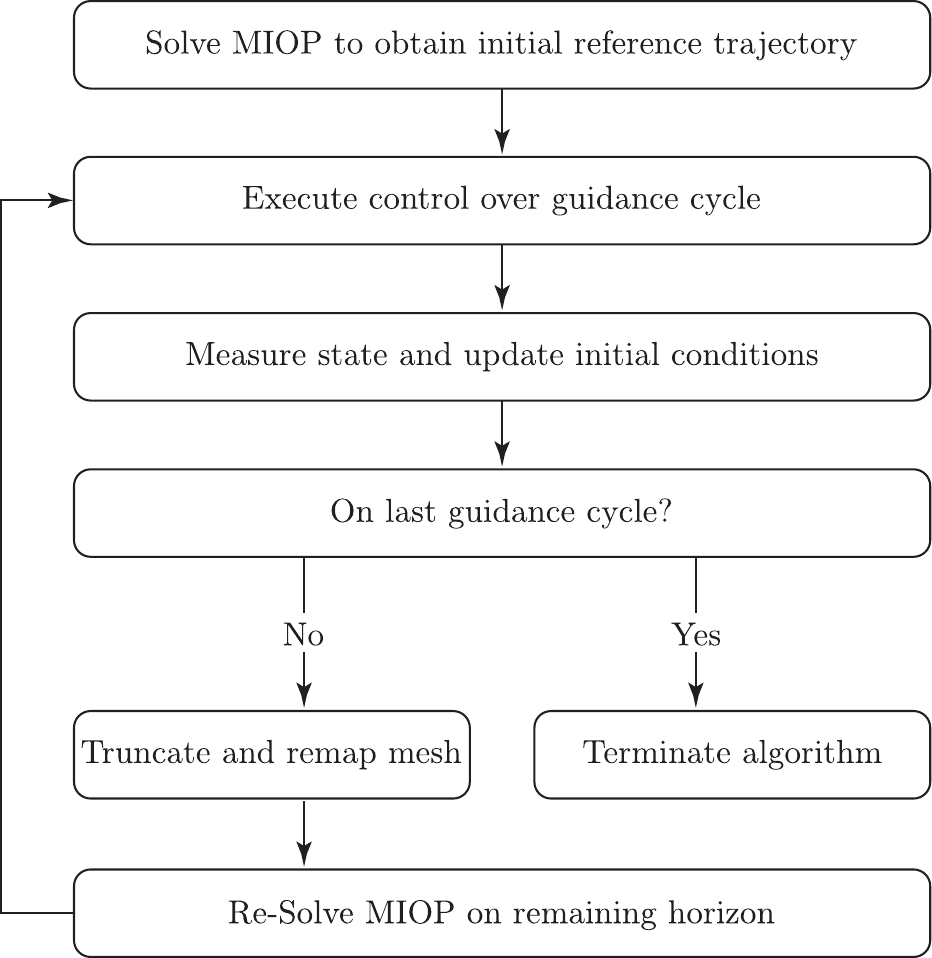}
\caption{Algorithmic flow diagram of the MIOG method.\label{fig:flowd-diagram}}
\end{figure}

Next, let $T_d^{(s)}$ be the time delay on guidance cycle $s$ and also be the computational time required to solve the MIOP on the remaining horizon $[t_0^{(s)},t_f]$.  The control that results from solving the MIOP on the remaining horizon $t\in[t_0^{(s)},t_f]$ is then used in the actual system (again, the system modeled in Eq.~\eqref{eq:actual-dynamics}) on $t\in[t_0^{(s)}+T_d^{(s)},t_e^{(s)}+T_d^{(s+1)}]$, with the initial state $\tilde{\m{x}}(t_e^{(s-1)}+T_d^{(s)})$. The process of solving the MIOP on the remaining horizon and using the resulting control in the actual system on the ensuing guidance cycle is repeated until the remaining horizon $t\in[t_0^{(s)},t_f]$ is shorter than the guidance cycle duration $\Delta T$.  Figure~\ref{fig:blockdiagram} provides a control system block diagram while Fig.~\ref{fig:flowd-diagram} provides the algorithmic flow of the MIOG method.

\subsection{Mesh Truncation and Remapping}\label{sect:mesh-remapping}

The margin-inducing guidance and control method of this paper requires that the LGR NLP associated with the MIOP be re-solved at the start of each guidance cycle.  Before the LGR NLP is solved at the start of  guidance cycle $s+1$, the mesh associated with guidance cycle $s$
 is truncated and remapped using the method developed in Ref.~\cite{Dennis2019}. This section provides a summary of the mesh truncation and remapping method.   Further details of the mesh truncation and remapping method can be found in Ref.~\cite{Dennis2019}. 

Recall that $t_0^{(s)}$ and $t_e^{(s)}$ are the start and end of guidance cycle $s$. Then, at the end of guidance cycle $s$, the time interval $t\in[t_0^{(s)},t_f]$ consists of an {\em expired} horizon, $t \in [t_0^{(s)},t_e^{(s)}]$, 
and an {\em unexpired} horizon, $t \in [t_e^{(s)},t_f]$. As a result, the collocation points associated with times less than $t_e^{(s)}$ are part of the expired horizon, while the collocation points corresponding to times greater than or equal to $t_e^{(s)}$ are part of the unexpired horizon. In the mesh truncation and remapping procedure of Ref.~\cite{Dennis2019}, the collocation points corresponding to the expired horizon are deleted from the mesh, while the mesh points corresponding to the unexpired horizon are retained.  Furthermore, the first collocation point on the unexpired horizon is moved to $t_e^{(s)}$ such that $t_0^{(s+1)}=t_e^{(s)}$ because this first collocation point on the unexpired horizon is greater than or equal to $t_e^{(s)}$.  After relocating the first collocation point on the unexpired horizon, the unexpired horizon $[t_e^{(s)},t_f]$ is transformed to  $\tau\in[-1,+1]$ for use with the LGR collocation method (see Section~\ref{sect:lgr-nlp}) and this truncated mesh on $\tau\in[-1,+1]$ is then a well aligned starting mesh for solving the MIOP on guidance cycle $s+1$ \cite{Dennis2019}. Finally, the initial state for solving the LGR NLP on guidance cycle $s+1$ is set to be the actual state $\tilde{\m{x}}(t_e^{(s)})$, that is, $\m{x}(t_0^{(s+1)})=\tilde{\m{x}}(t_e^{(s)})$.

\begin{figure}[h]
\centering
\includegraphics[scale=0.75]{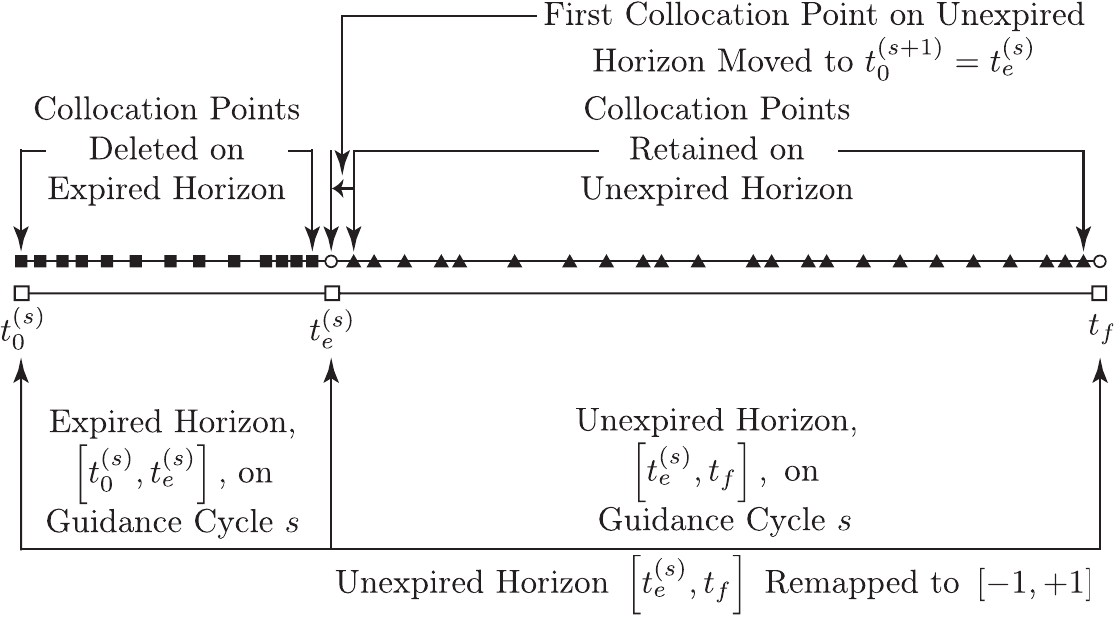}
\caption{Schematic of Mesh Truncation and Remapping Method.\label{fig:mesh-truncation}}
\end{figure}

The process of deleting the expired portion of the mesh while retaining the unexpired portion of the mesh is repeated at the start of each guidance cycle.  As a result, the mesh size progressively shrinks as the number of guidance cycles increases. A schematic of the mesh truncation and remapping procedure is provided in Fig.~\ref{fig:mesh-truncation}. Further details on the mesh truncation and remapping procedure can be found in Ref.~\cite{Dennis2019}. 

\subsection{Infeasible or Time-Consuming NLP}\label{sect:Infeas-nlp}

While the MIOP does not contain state path constraints, it is  possible on any guidance cycle that the LGR NLP described in Section \ref{sect:lgr-nlp} may not be feasible or may be excessively time consuming, where the definition of excessively time consuming is that the NLP takes longer than $\Delta T$ time units to solve (where it is recalled that $\Delta T$ is the guidance cycle duration). In this paper the approach taken to safeguard against these two possibilities is to continue to use the most recently computed control.  That is, suppose that a solution has been obtained to the LGR NLP given in Section~\ref{sect:lgr-nlp} at guidance update time $t_0^{(s)}$.  Furthermore, suppose that for $n\geq 1$ consecutive guidance cycles with guidance update times
$(t_0^{(s+1)},\ldots,t_0^{(s+n)})$ either the NLP solver is unable to find a solution or the NLP solver takes greater than $\Delta T$ time units to find a solution, and that a solution to the LGR NLP is again obtained at guidance update time $t_0^{(s+n+1)}$ within $\Delta T$ time units.  Then the control obtained on $[t_0^{(s)},t_f]$ is used in the actual system until the guidance update time $t_0^{(s+n+1)}$, and a new control is obtained at $t_0^{(s+n+1)}$. 

\section{Formulation of MIOP for RLV Entry}\label{sect:rlv-ocp}

This section formulates the margin-inducing optimal control problem (MIOP) for the reusable launch vehicle (RLV) entry optimal control problem which is used to demonstrate the margin-inducing optimal guidance (MIOG) method.  Specifically, the original path-constrained optimal control problem under consideration includes state path constraints on heating rate and dynamic pressure. The objective of this RLV entry problem is to improve crossrange while attenuating phugoid oscillations and maintaining control authority during entry.  The constrained optimal control problem is reformulated into a MIOP of the form described in Section \ref{sect:miop}. 

\subsection{RLV Dynamics}\label{sect:rlv-equations-of-motion}

The reusable launch vehicle (RLV) dynamics correspond to a vehicle in motion over a spherical nonrotating planet. These dynamics are divided into the state dynamics
\begin{equation}\label{eq:rlv-eom}
  \begin{array}{lcl}
    \dot{r} & = & \displaystyle v\sin\gamma, \vspace{3pt} \\
    \dot{\theta} & = & \displaystyle \frac{v\cos\gamma\sin\psi}{r\cos\phi},\vspace{3pt} \\
    \dot{\phi} & =& \displaystyle \frac{v\cos\gamma\cos\psi}{r},\vspace{3pt} \\
    \dot{v} & = & \displaystyle -D-g\sin\gamma , \vspace{3pt} \\
    \dot{\gamma}& =& \displaystyle \frac{1}{v}\left[L\cos\sigma-\left(g-\frac{v^2}{r}\right)\cos\gamma\right] , \vspace{3pt} \\
    \dot{\psi}& =& \displaystyle \frac{L\sin\sigma}{v\cos\gamma}+\frac{v\cos\gamma\sin\psi\tan\phi}{r}, 
  \end{array}
\end{equation}
and the control (actuator) dynamics
\begin{equation}\label{eq:rlv-actuator-dynamics}
\begin{array}{lcl}
    \dot{\alpha} & = & p_\alpha, \\
    \dot{\sigma} & = & p_\sigma,
\end{array}
\end{equation}
where $r$ is the geocentric radius, $v$ is the speed, $\theta$ is the longitude, $\phi$ is the latitude, $\gamma$ is the flight path angle,  $\psi$ is the azimuth (measured clockwise from due north), $g = \mu/r^2$ is the gravitational acceleration, $\mu$ is the gravitational parameter, $q=\rho v^2/2$ is the dynamic pressure, $\rho=\rho_0\exp(-h/H)$ is the atmospheric density, $h = r - R_e$ is the altitude, $R_e$ is the radius of the Earth, $\rho_0$ is the surface level atmospheric density, $H$ is the density scale height,  $D=qSC_D/m$ is the drag specific force, $L=qSC_L/m$ is the lift specific force, $C_D=C_{D0}+C_{D1}\alpha+C_{D2}\alpha^2$ is the coefficient of drag, $C_L=C_{L0}+C_{L1}\alpha$ is the lift coefficient, $\alpha$ is the angle-of-attack, $\sigma$ is the bank angle, $p_\alpha$ is the angle of attack rate, and $p_\sigma$ is the bank angle rate.  It is noted that $(p_\alpha,p_\sigma)$ forms the pseudo-control for this problem.  Finally, the reference aerodynamic coefficients and physical parameters used in this study are given in Table \ref{tab:physical-constants-Earth}.

\begin{table}[ht]
\centering
\caption{Nominal physical and aerodynamic constants \label{tab:physical-constants-Earth}}
\renewcommand{\baselinestretch}{1.25}\normalsize\normalfont
\begin{tabular}{ccc}\hline\hline
Quantity & Value & Unit  \\\hline
  $R_e$ & $6.371\times10^6$ & m   \\
  $r_n$ & $1$ & m \\
  $\mu$ & $3.986\times10^{14}$ & $\textrm{m}^3$/$\textrm{s}^2$ \\
  $g_0$ & $9.80665$ & m/$s^2$ \\
  $S$ & $249.9$ & $\textrm{m}^2$ \\
   $K_q$ & $1.7415\times10^{-4}$ & kg${}^{1/2}$/m \\
  $C_{D0}$ & $0.0785$ & - \\
   $C_{D1}$ & $-0.3529 $ & - \\
    $C_{D2}$ & $2.0400$ & -  \\
  $C_{L0}$ & $-0.207$ & - \\
  $C_{L1}$ & $1.6756$ & - \\
  $m$ & $92079$ & kg  \\
  $\rho_0$ & $1.2256$ & kg/$\textrm{m}^3$\\
   $H$ & $7254$ &m  \\\hline\hline
\end{tabular}
\end{table}

\subsection{RLV Boundary Conditions}\label{sect:rlv-boundary-conditions}

The boundary conditions for the RLV entry problem are given as
\begin{equation}\label{eq:rlv-boundary-conditions}
\begin{array}{lclclcl}
    r(t_0) & = & r_0 & , & r(t_f) & = & r_f,    \\
    \theta(t_0) & = & \theta_0 & , & \theta(t_f) & = & \textrm{Free}, \\
    \phi(t_0) & = & \phi_0  & , & \phi(t_f) & = & \textrm{Free}, \\
    v(t_0) & = & v_0  & , & v(t_f) & = & v_f, \\
    \gamma(t_0) & = & \gamma_0 & , & \gamma(t_f) & = & \gamma_f, \\
    \psi(t_0) & = & \psi_0 & , & \psi(t_f) & = & \textrm{Free}, \\
    \alpha(t_0) & = & \textrm{Free} & , & \alpha(t_f) & = & \textrm{Free}, \\
    \sigma(t_0) & = & \textrm{Free} & , & \sigma(t_f) & = & \textrm{Free}.
    \end{array}
\end{equation}
Finally, it is noted that, while the initial angle of attack and bank angle, $(\alpha(t_0),\sigma(t_0))$, are free, on any guidance cycle $s>0$ the angle of attack and bank angle are continuous, that is
\begin{equation}
\left[
  \begin{array}{c}
    \alpha(t_0^{(s+1)}) \\
    \sigma(t_0^{(s+1)})
  \end{array}
\right] = 
\left[
  \begin{array}{c}
    \alpha(t_e^{(s)}) \\
    \sigma(t_e^{(s)})
  \end{array}
\right], \quad (s=1,\ldots,S-1).
\end{equation}
The initial conditions and terminal conditions for the RLV entry are provided in Table \ref{tab:bcs-Earth}.

\begin{table}[ht]
\centering
\caption{Boundary conditions \label{tab:bcs-Earth}}
\renewcommand{\baselinestretch}{1.25}\normalsize\normalfont
\begin{tabular}{ccc}\hline\hline
  Quantity& Initial & Terminal\\\hline
  $r$ & $6450$~km  & $6396$~km \\
  $\theta$ & $0$~deg  & Free \\
  $\phi$ & $0$~deg & Free \\
   $v$ & $7.80288$~km/s  & $0.762$~km/s \\
  $\gamma$ & $-1$~deg  & $-5$~deg \\
  $\psi$ & $90$~deg  & Free \\ \hline\hline
\end{tabular}
\end{table}

\subsection{RLV Entry Path Constraints}\label{sect:rlv-path-constraints}

The constraints in the RLV entry problem are divided into control constraints and state constraints, where it is noted that the state constraints are those for which margin is to be induced during guidance.  

The control, pseudo-control, and state path constraints during entry are given, respectively, as
\begin{equation}\label{eq:rlv-control-bounds}
    \begin{array}{lcccl}
        \alpha_{\min} & \leq & \alpha & \leq & \alpha_{\max}, \\
        \sigma_{\min} & \leq & \sigma & \leq & \sigma_{\max},
    \end{array}
\end{equation}
\begin{equation}\label{eq:rlv-pseudo-control-bounds}
    \begin{array}{lcccl}
        p_{\alpha,\min} & \leq & p_\alpha \leq & p_{\alpha,\max}, \\
        p_{\sigma,\min} & \leq & p_\sigma \leq & p_{\sigma,\max},
    \end{array}
\end{equation}
and
\begin{equation}
\begin{array}{lcl}
    \dot{Q} & \leq & \dot{Q}_{\max}, \\
        q   & \leq & q_{\max},
\end{array}
\end{equation}
where $\dot{Q}=K_q \sqrt{\rho/r_n} v^3$ is the heating rate \cite{sutton1971general}, $K_q$ is a heating rate constant, and $r_n$ is the nose radius of the entry vehicle. The constraint limits are provided in Table~\ref{tab:constraint-limits}.

\begin{table}[ht]
\centering
\caption{Constraint limits \label{tab:constraint-limits}}
\renewcommand{\baselinestretch}{1.25}\normalsize\normalfont
\begin{tabular}{ccc}\hline\hline
  Quantity & Minimum & Maximum\\\hline
  $\alpha$~(deg) &$-90$  & $90$ \\
  $\sigma$~(deg) & $-90$  & $+90$ \\
  $p_{\alpha}$~(deg/s) &  $-0.5$ & $0.5$ \\
  $p_{\sigma}$~(deg/s) &  $-5$ & $5$ \\
  $\dot{Q}$~(MW/m${}^2$) & $0$  & $0.85$ \\
  $q$~(kPa) & $0$  & $16.375$ \\\hline\hline
\end{tabular}
\end{table}

\subsection{RLV Entry Objective Functional}\label{sect:rlv-objective}

The objective functional, $\C{J}$, for the constrained RLV entry problem is given as
\begin{equation}\label{eq:rlv-objective}
\C{J} = -c_1\phi(t_f)+ \int_{t_0}^{t_f} \left\{ c_2\left[\left(\frac{p_{\alpha}}{p_{\alpha,\max}}\right)^2+\left(\frac{p_{\sigma}}{p_{\sigma,\max}}\right)^2\right]+ c_3\left[\coth{\left(-\frac{k}{2}\right)}\tanh\left(\frac{-k\sin\gamma}{2}\right)\right]^2 \right\} dt.
\end{equation}
The objective functional in Eq.~\eqref{eq:rlv-objective} is designed such that the terms lead to a well-behaved guidance trajectory.   The first term in Eq.~\eqref{eq:rlv-objective} is included to maximize the terminal cross-range and is the original objective functional shown in Ref.~\cite{Betts2020}. The second term in Eq.~\eqref{eq:rlv-objective} is included to maintain a smooth control \cite{miller2022end}. The third term in the objective functional is designed to reduce phugoid oscillations during atmospheric entry \cite{miller2022end}. The term intended to prevent phugoid oscillations is designed such that when $\gamma = 0$, the term is zero. This term monotonically increases to unity as the flight path angle increases or decreases from $\gamma = 0$ to $\gamma = \pm\pi/2$, where $k$ is a constant that determines how quickly the term increases.

\subsection{Margin-Inducing Objective Functional}\label{sect:rlv-margin-inducing-objective}

As described in Section \ref{sect:miop}, the margin-inducing objective functional includes terms intended to increase the margin between the path constraint functions and their limits.  In the case of the RLV entry optimal control problem under consideration in this section, these additional terms involve the heating rate and dynamic pressure path constraint functions.  Consequently, in terms of Eq.~\eqref{eq:objective-miop} which in turn includes terms of the form given in Eq.~\eqref{eq:max-approximation-continuous}, the RLV entry margin-inducing objective functional 
is given as
\begin{equation}\label{eq:rlv-margin-inducing-objective}
\C{J}_{\textrm{a}}= \C{J}+c_{\dot{Q}}\C{E}_{\beta}(e_{\dot{Q}}) + c_q\C{E}_{\beta}(e_q),
\end{equation}
where $e_{\dot{Q}}=\dot{Q}/\dot{Q}_{\max}$ and $e_q=q/q_{\max}$ are the scaled state path constraint functions.  

\subsection{RLV Entry Margin-Inducing Optimal Control Problem}\label{subsect:rlv-miop}

The margin-inducing RLV entry optimal control problem is then stated as follows: minimize the objective functional given in Eq.~\eqref{eq:rlv-margin-inducing-objective} subject to the differential equation constraints given in Eq.~\eqref{eq:rlv-eom}, the control rate (actuator) dynamics of Eq.~\eqref{eq:rlv-actuator-dynamics}, the boundary conditions given in Eq.~\eqref{eq:rlv-boundary-conditions}, the control constraints given in Eq.~\eqref{eq:rlv-control-bounds}, and the pseudo-control (actuator) constraints given in Eq.~\eqref{eq:rlv-pseudo-control-bounds}.  The objective functional parameters used in the margin-inducing RLV entry optimal control problem are given in Table~\ref{tab:obj-functional-parameters}.

\begin{table}[ht]
\centering
\caption{Objective functional parameters  \label{tab:obj-functional-parameters}}
\renewcommand{\baselinestretch}{1.25}\normalsize\normalfont
\begin{tabular}{cc}\hline\hline
  Parameter & Value \\\hline
  $c_1$ & $0.8$ \\
  $c_2$ & $0.01$ \\
  $c_3$ & $0.5$ \\
   $c_{\dot{Q}}$ & $0.2$ \\
     $c_q$ & $2$ \\
   $k$ & $3$ \\
 $\beta$ & $10$ \\       
  \hline\hline
\end{tabular}
\end{table}

\section{Application of Margin-Inducing Guidance and Control to RLV Entry}\label{sect:rlv-results-miog}


This section presents the results of employing both the MIOG method of this paper and the computational optimal guidance (COG) method of Ref.~\cite{Dennis2019} on the RLV entry problem described in Section \ref{sect:rlv-ocp} using a guidance cycle duration of $\Delta T = 15$~s, where it is noted that the COG method solves the original constrained optimal control problem (COP) described in Section \ref{sect:constrained-ocp}.  It was assumed in this study that uncertainties exist in the parameters $(m,\rho_0,C_D,C_L)$ and the initial state $\m{x}(t_0)$.  Moreover, as seen in Table~\ref{tab:modeling-uncertainties}, uncertainties in $(m,\rho_0,C_D,C_L)$ were uniformly distributed while uncertainties in $\m{x}(t_0)$ were normally distributed, where the distributions of these uncertainties are shown in Table~\ref{tab:modeling-uncertainties} and the means of these uncertainties are given in Tables~\ref{tab:physical-constants-Earth} and \ref{tab:bcs-Earth}.  The performance of the COG method was evaluated using $50$ Monte Carlo trials, 
where $50$ trials were used for COG because infeasibility was readily apparent within this sample size.  The performance of the MIOG method was then evaluated using $1000$ Monte Carlo trials.  

The results shown in this paper were obtained using the $hp$-adaptive LGR collocation optimal control software $\mathbb{GPOPS-II}$ \cite{Garg2010,Patterson2014,Patterson2015,Garg2011a,Garg2011b} with the NLP solvers IPOPT \cite{Biegler2008} (with the linear solver MA57 \cite{Duff2004MA57}) and SPRNLP \cite{Betts1994,Betts-Rao-2026-SOS-NLP}.  All first and second derivatives required by the NLP solvers were obtained using the algorithmic differentiation software ADiGator \cite{Weinstein2017}.  Furthermore, any mesh refinement required by $\mathbb{GPOPS-II}$ was performed using the method described in Ref.~\cite{Liu2018}. When computing the  MIOP or COP reference solutions, the mesh refinement error tolerance was $10^{-5}$ with a maximum number of ten (10) mesh refinement iterations.  When implementing MIOG or COG for solving the guidance and control problem, a mesh refinement error tolerance of $10^{-5}$ was used with a maximum of two mesh refinement iterations.  Furthermore, all simulations of the actual dynamics were performed using the \textsf{MATLAB} differential equation solver \textsf{ode113} with a relative error tolerance of $10^{-8}$ and a cubic spline interpolant of the LGR NLP control. Finally, all results presented in this paper were obtained using MATLAB R2026a (Build 
26.1.0.3312084) on an Apple MacBook Pro with an M5 Pro processor and 48 GB RAM running MacOS Tahoe Version 26.6.2 (Build 25G83).  

\begin{table}[h]
    \centering
    \caption{Modeling Uncertainties}
    \label{tab:modeling-uncertainties}
    
    \begin{subtable}[t]{\textwidth}
        \centering

\begin{tabular}{cc}\hline\hline
  Parameter & Value \\\hline
  $m$ & $[-5\%,+5\%]$ \\
  $\rho_0$ & $[-10\%,+10\%]$ \\
  $C_D$ & $[-10\%,+10\%]$ \\
   $C_L$ & $[-10\%,+10\%]$ \\
  \hline\hline
\end{tabular}

        \caption{Parameter Uncertainties}
        \label{tab:uniform-distributions}
    \end{subtable}
\\
    \begin{subtable}[t]{\textwidth}
        \centering

\begin{tabular}{cc}\hline\hline
  Initial State & 3-$\sigma$ \\\hline
  $r_0$~(m) & $100$ \\
  $\theta_0$~(deg) & $0.05$ \\
  $\phi_0$~(deg) & $0.05$ \\
   $v_0$~(m/s) & $20$ \\
   $\gamma_0$~(deg) & $0.05$ \\
   $\psi_0$~(deg) & $0.05$ \\
  \hline\hline
\end{tabular}

        \caption{Initial State Uncertainties}
        \label{tab:IC perts}
    \end{subtable}

\end{table}

\subsection{Results Using COG Method of Ref.~\cite{Dennis2019}}

Before proceeding to the performance of the MIOG method, this section provides results using the COG method of Ref.~\cite{Dennis2019} using $50$ Monte Carlo trials.  This reduced set of Monte Carlo trials was sufficient to demonstrate the issues with the COG method for a problem with state path constraints.

\begin{figure}[htbp]
    \centering

    \begin{minipage}[t]{0.48\textwidth}
        \begin{subfigure}[t]{\linewidth}
            \centering
            \includegraphics[width=\linewidth]{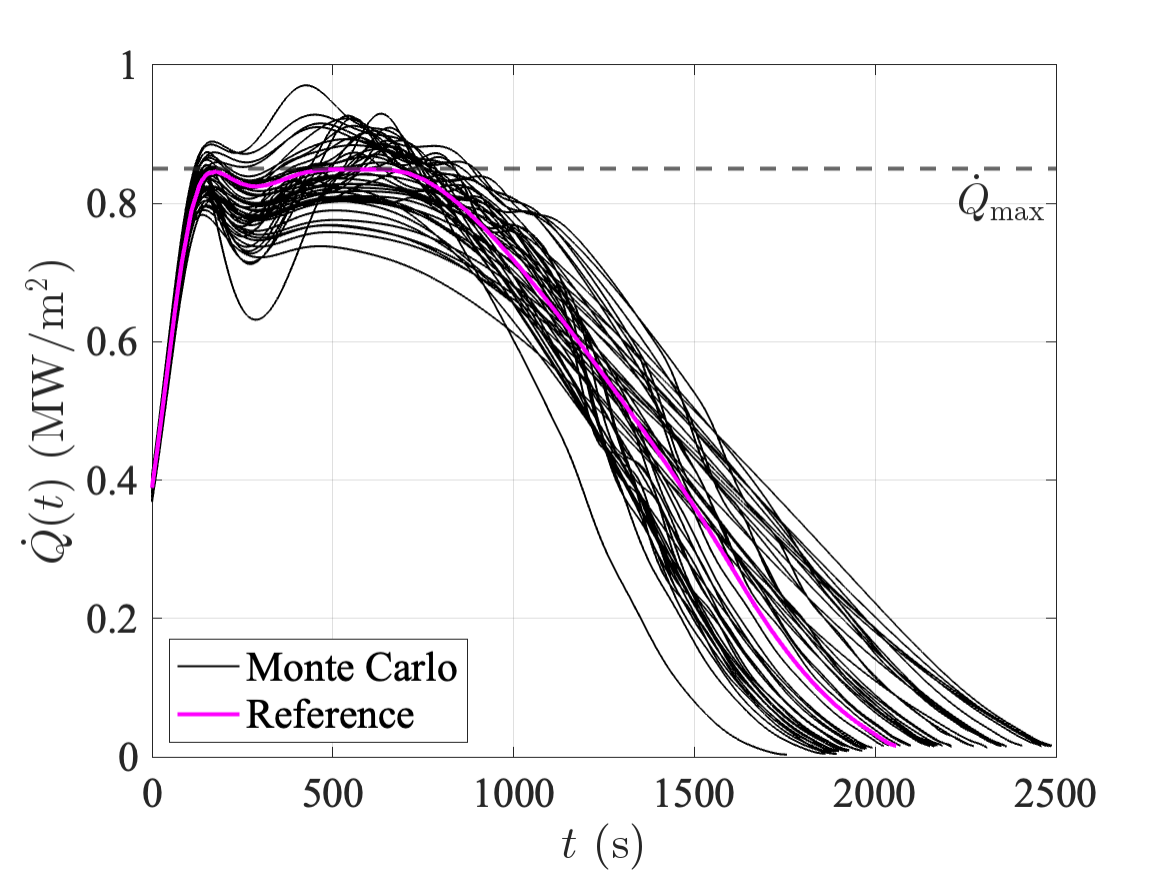}
            \caption{ $\dot{Q}(t)$ vs  $t$ using IPOPT.}
            \label{fig:Qdot-IPOPT-COG}
        \end{subfigure}

        \vspace{1em}

        \begin{subfigure}[t]{\linewidth}
            \centering
            \includegraphics[width=\linewidth]{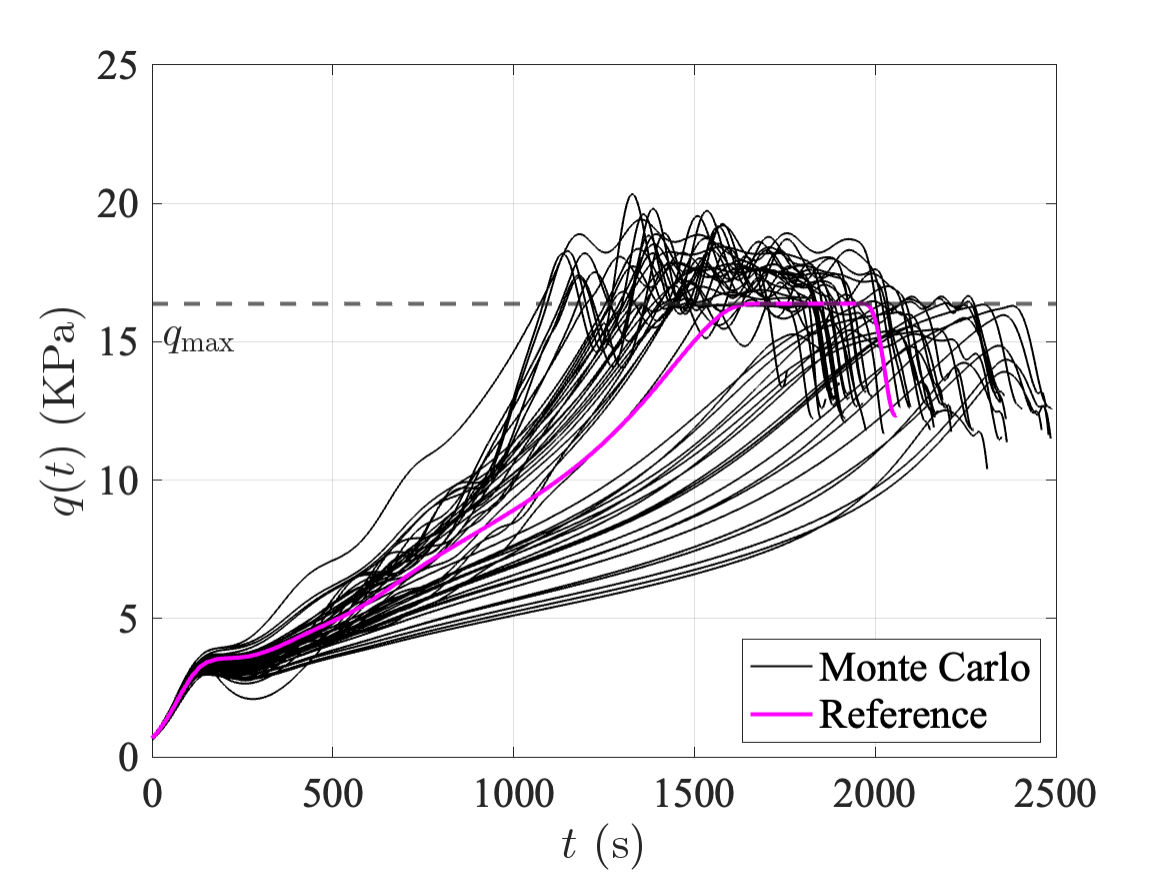}
            \caption{ $q (t) $ vs $t$ using IPOPT.}
            \label{fig:q-IPOPT-COG}
        \end{subfigure}

    \end{minipage}
    \hfill
    \begin{minipage}[t]{0.48\textwidth}
        \begin{subfigure}[t]{\linewidth}
            \centering
            \includegraphics[width=\linewidth]{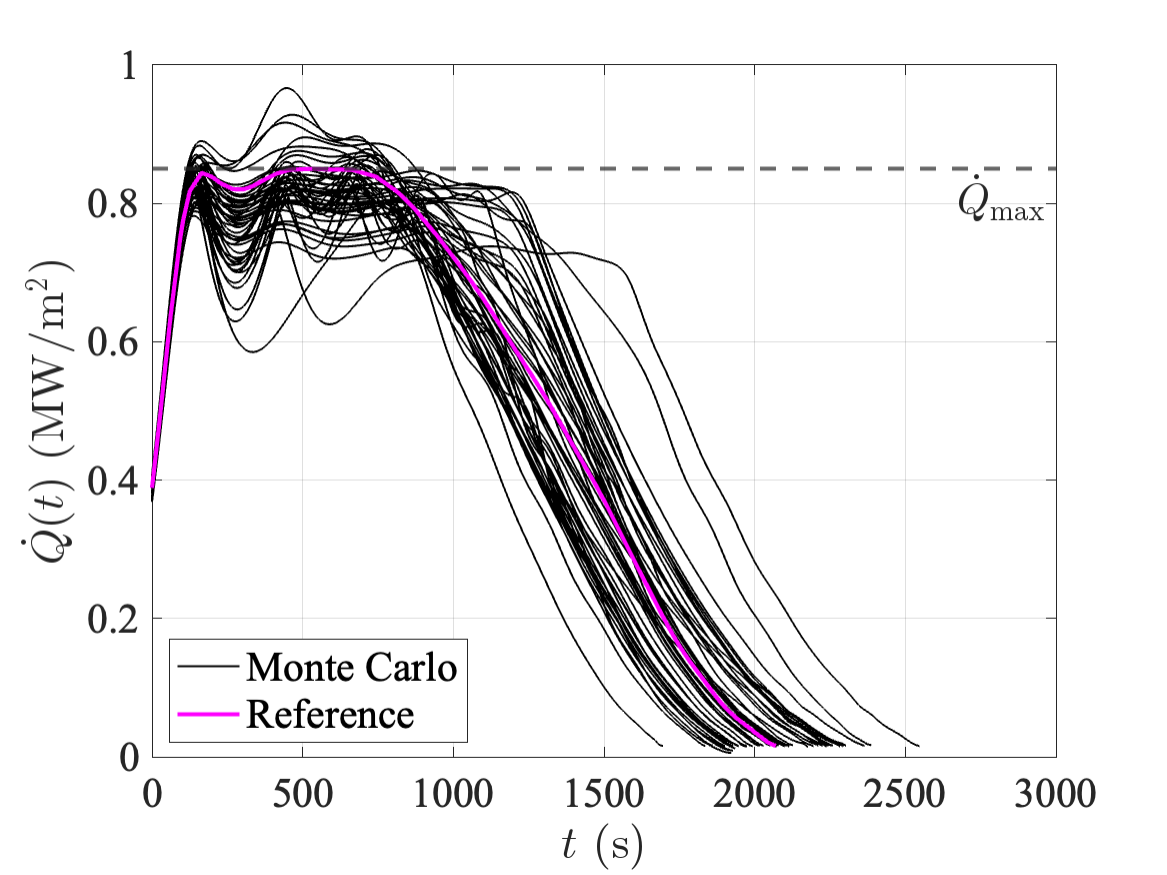}
            \caption{ $\dot{Q}(t)$ vs $t$ using SPRNLP.}
            \label{fig:Qdot-SOS-COG}
        \end{subfigure}

        \vspace{1em}

        \begin{subfigure}[t]{\linewidth}
            \centering
            \includegraphics[width=\linewidth]{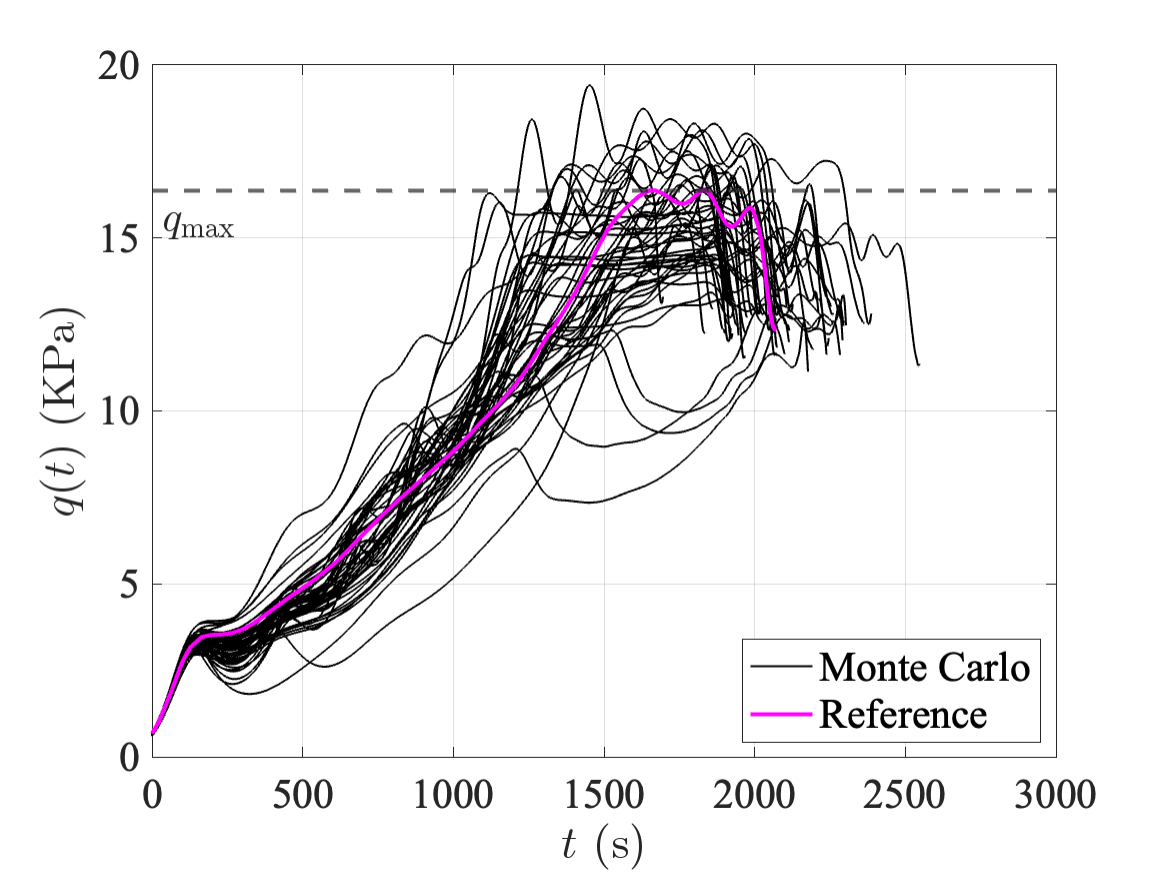}
            \caption{ $q (t) $ vs $t$ using SPRNLP.}
            \label{fig:q-SOS-COG}
        \end{subfigure}

    \end{minipage}

    \caption{Heating rate and dynamic pressure obtained using the COG method with NLP solvers IPOPT [(a)--(b)] and SPRNLP [(c)--(d)].}
    \label{fig:heating-rate-dyn-pressure-COG}
\end{figure}

Figure~\ref{fig:heating-rate-dyn-pressure-COG} shows the results of applying the COG method using both IPOPT and SPRNLP.  It is seen that the path constraints are violated frequently using the COG method, and the amounts of the violations are quite significant.  A detailed breakdown of the statistics of these constraint violations is shown in 
Table~\ref{tab:max-cons-Cog} for the maximum heating rate and dynamic pressure. As seen in Table~\ref{tab:max-cons-Cog}, the mean values of $\max(\dot{Q}),\max(q)$ are $(1.41\%,3.81\%)$ above the constraint limit using IPOPT and $(0.88\%,4.20\%)$ below the limit using SPRNLP.  In addition, the maximum $\max(\dot{Q}),\max(q)$ are $(14.19\%,24.17\%)$ above the constraint limit using IPOPT and $(13.73\%,18.56\%)$ above the constraint limit using SPRNLP.  Moreover, as shown in Table~\ref{tab:feas-results-COG}, these constraint violations lead to a significant number of guidance cycles where a solution to the COP cannot be obtained. In the worst cases, an NLP solution was not obtained on $93$ and $83$ guidance cycles using the COG method with IPOPT and SPRNLP, respectively. Finally, it was found that $62\%$ of the IPOPT and $44\%$ of the SPRNLP guidance trajectories violate the state path constraints during guidance.  The results provided in this section demonstrate that the COG method is not well-suited to problems with state path constraints.  

\begin{table}[htbp]
\centering
\caption{Maximum heating rate and dynamic pressure results using the COG method (margin amounts in parentheses).}
\label{tab:max-cons-Cog}
\centering

\begin{tabular}{lcccc}
\hline\hline
& \multicolumn{2}{c}{IPOPT}
& \multicolumn{2}{c}{SPRNLP} \\
\cline{2-3}\cline{4-5}
Statistic
& $\max(\dot{Q})$ & $\max(q)$
& $\max(\dot{Q})$ & $\max(q)$ \\
& (MW/m$^2$) & (kPa)
& (MW/m$^2$) & (kPa) \\
\hline
Reference
& $0.85\;(0)$ & $16.375\;(0)$
& $0.85\;(0)$ & $16.375\;(0)$ \\
Mean
& $0.862\;(-0.012)$ & $16.999\;(-0.624)$
& $0.843\;(0.0074)$ & $15.687\;(0.688)$ \\
Max
& $0.971\;(-0.121)$ & $20.333\;(-3.958)$
& $0.967\;(-0.117)$ & $19.414\;(-3.039)$ \\
Std. Dev.
& $0.0447\;(0.0447)$ & $1.694\;(1.694)$
& $0.0353\;(0.0353)$ & $1.559\;(1.559)$ \\
\hline\hline
\end{tabular}

\end{table}

\begin{table}[htbp]
\centering
\caption{COG guidance cycles where an NLP solution was not obtained.}
\label{tab:feas-results-COG}

\begin{tabular}{lcc}\hline\hline
Statistic
& IPOPT & SPRNLP\\
\hline
Mean
& $43.1$ & $23.2$ \\
Max
& $93$ & $83$ \\
Std. Dev.
& $37.9$ & $26.3$ \\
\hline\hline
\end{tabular}

\end{table}

\subsection{Results Using MIOG Method}\label{sect:rlv-margin-inducing-results}

This section provides the key results obtained using the MIOG method.  As stated above, the results using the MIOG method were generated with $1000$ Monte Carlo trials. The results include those pertaining to margin-induction, behavior of the guidance trajectories, computational performance, and terminal state accuracy. 

\subsubsection{Margin-Inducing Feature of MIOG Method}

\begin{figure}[htbp]
    \centering

    \begin{minipage}[t]{0.48\textwidth}
        \begin{subfigure}[t]{\linewidth}
            \centering
            \includegraphics[width=\linewidth]{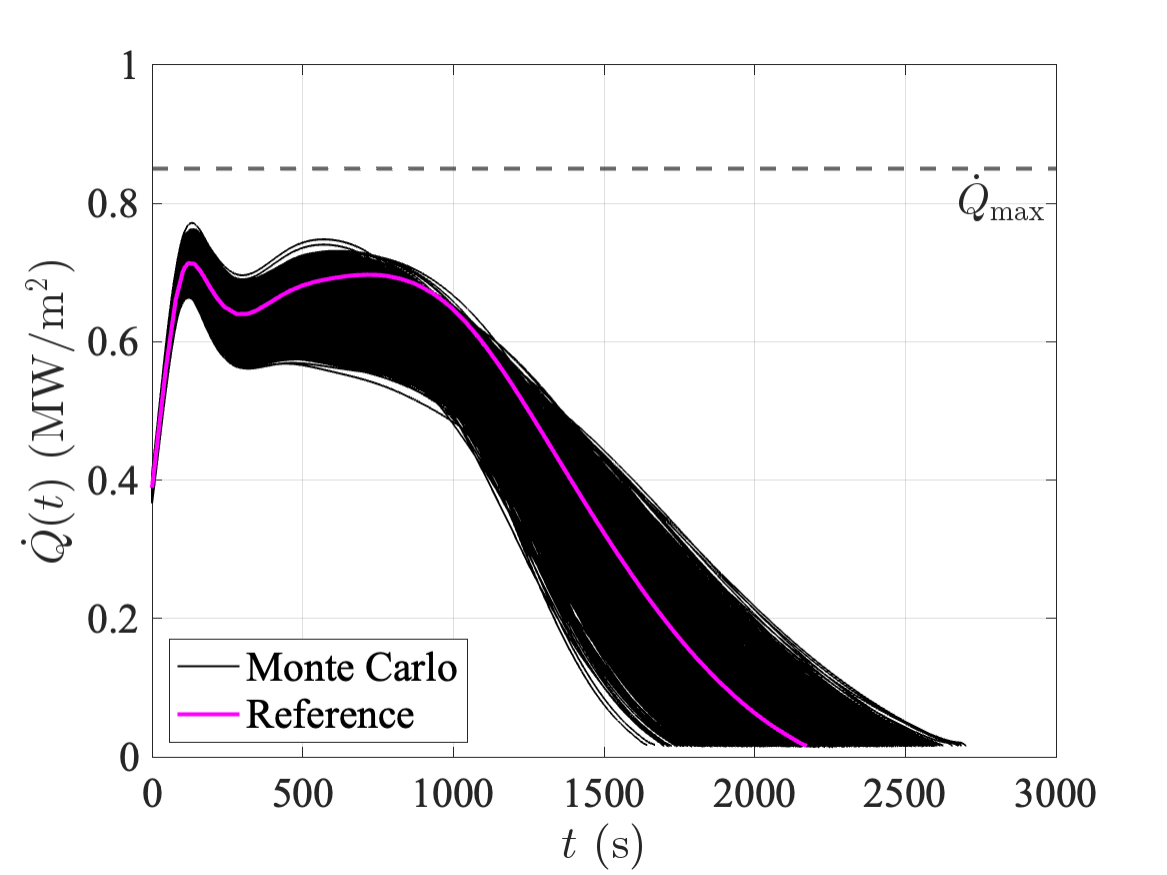}
            \caption{ $\dot{Q}(t)$ vs.~$t$ using IPOPT.}
            \label{fig:Qdot-IPOPT}
        \end{subfigure}

        \vspace{1em}

        \begin{subfigure}[t]{\linewidth}
            \centering
            \includegraphics[width=\linewidth]{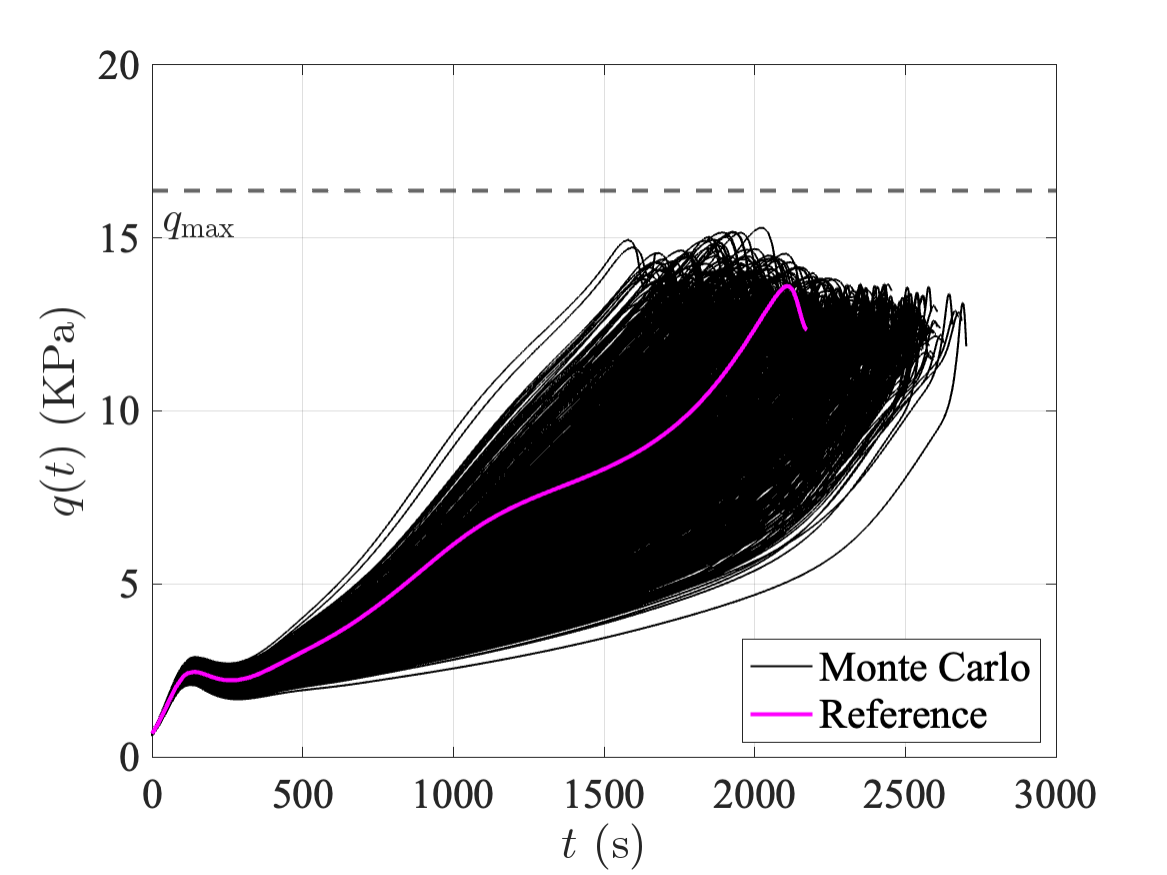}
            \caption{ $q(t)$ vs.~$t$ using IPOPT.}
            \label{fig:q-IPOPT}
        \end{subfigure}
    \end{minipage}
    \hfill
    \begin{minipage}[t]{0.48\textwidth}
        \begin{subfigure}[t]{\linewidth}
            \centering
            \includegraphics[width=\linewidth]{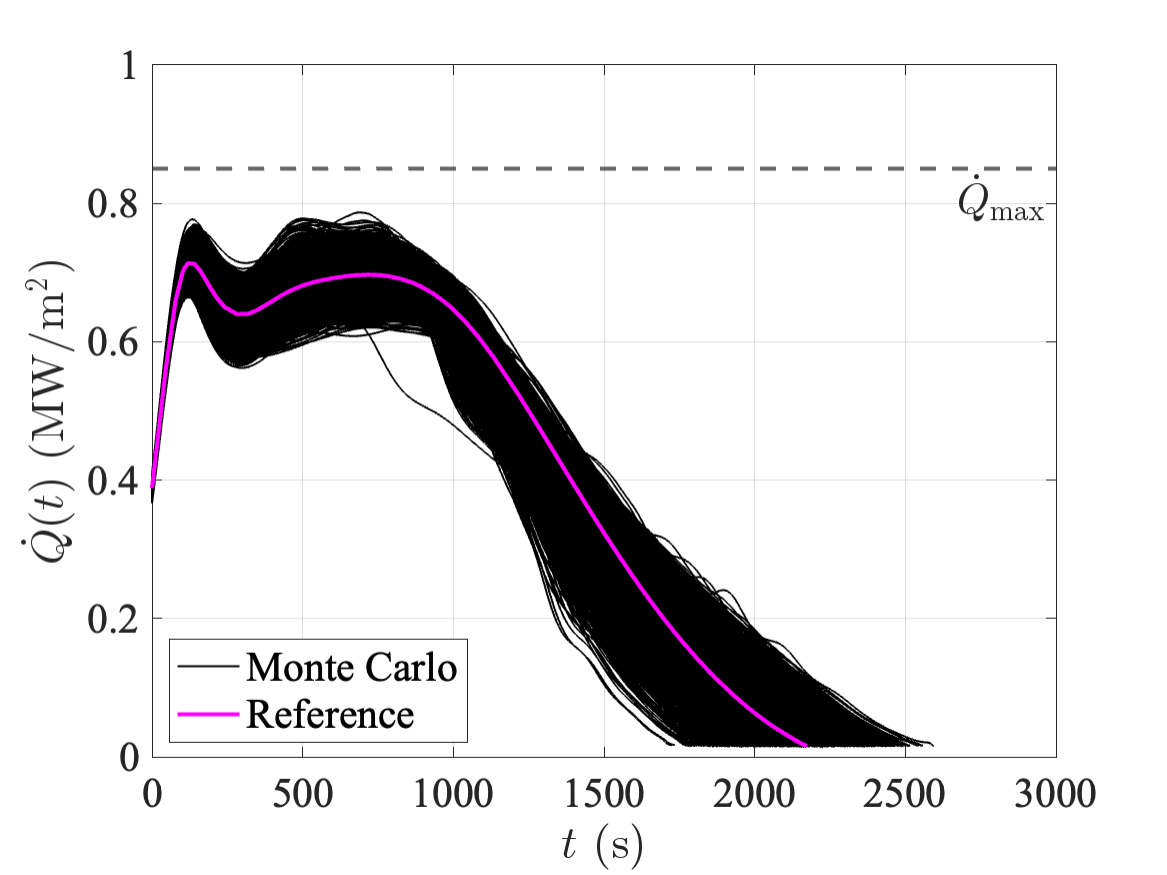}
            \caption{ $\dot{Q}(t)$ vs.~$t$ using SPRNLP.}
            \label{fig:Qdot-SOS}
        \end{subfigure}

        \vspace{1em}

        \begin{subfigure}[t]{\linewidth}
            \centering
            \includegraphics[width=\linewidth]{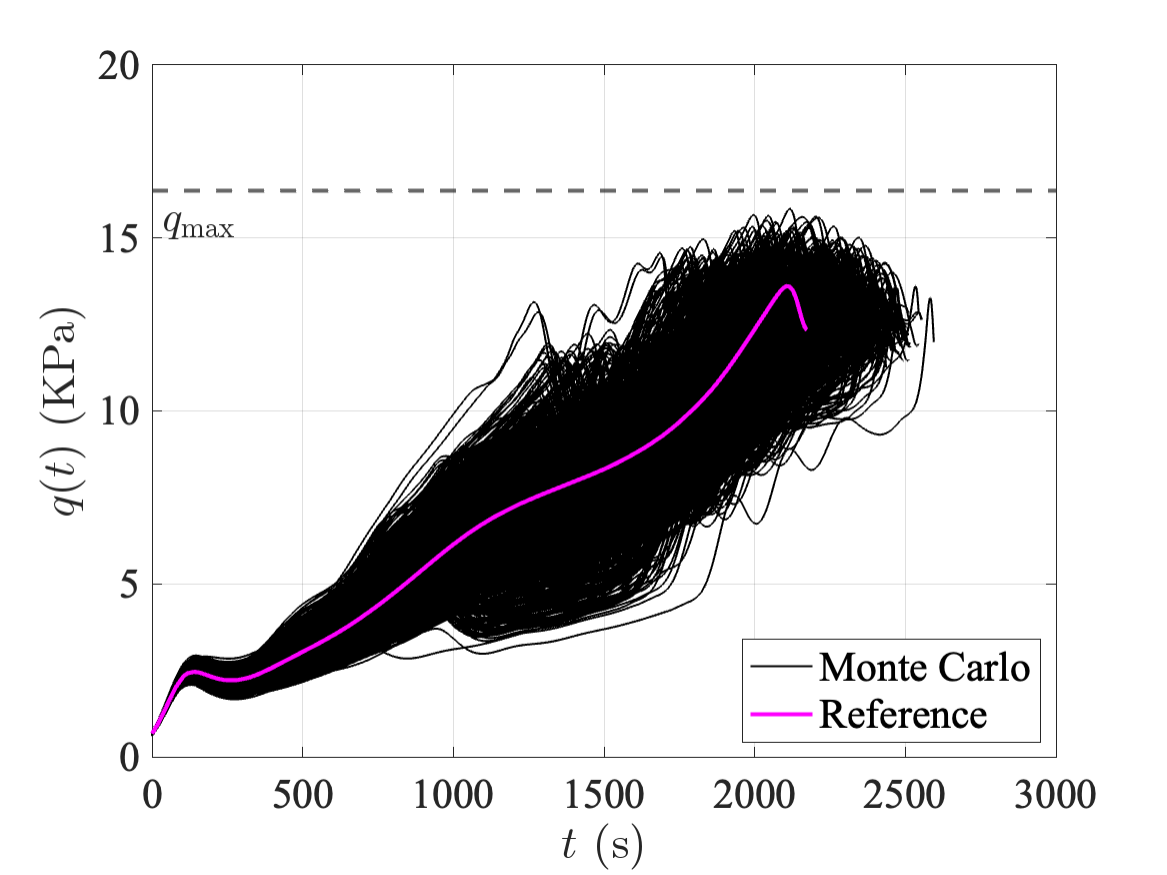}
            \caption{ $q(t) $ vs.~$t$ using SPRNLP.}
            \label{fig:q-SOS}
        \end{subfigure}       
    \end{minipage}

    \caption{Heating rate and dynamic pressure obtained using the MIOG method with NLP solvers IPOPT [(a)--(b)] and SPRNLP [(c)--(d)].}
    \label{fig:heating-rate-dyn-pressure-MIOG}
\end{figure}

The first key feature of the MIOG method is that the method induces margin in the state path constraints during guidance. Specifically, Figs.~\ref{fig:Qdot-IPOPT}--\ref{fig:q-SOS} show that the heating rate and dynamic pressure obtained from the MIOG method remain significantly below their respective constraint limits for all 1000 simulated trajectories. Table~\ref{tab:max-path-cons-MIOG} provides the statistical results for the maximum heating rate and dynamic pressure (and their margins) using IPOPT and SPRNLP.  First, it is seen that the reference heating rate and dynamic pressure margins are $16.05\%$ and $16.898\%$, respectively, for both NLP solvers. Second, the mean heating rate margins are $16.1\%$ (IPOPT) and $15.65\%$ (SPRNLP) while the mean dynamic pressure margins are $21.26\%$ (IPOPT) and $16.90 \%$ (SPRNLP). The minimum heating rate margins are $ 9.21\%$ (IPOPT) and  $7.43 \%$ (SPRNLP), and the minimum dynamic pressure margins are $6.60\%$ (IPOPT) and $ 3.25 \%$ (SPRNLP). These results demonstrate that the MIOG method is capable of maintaining margin in the state path constraints in the presence of a wide range of uncertainties.

In a manner consistent with inducing margin, the MIOG method produces solutions that remain feasible during guidance. In particular, Table~\ref{tab:feas-results} shows that an NLP solution to the MIOP was obtained on nearly every guidance cycle, and the NLP solver was unsuccessful on guidance cycles near the end of the guidance trajectories due to limited control authority on the short remaining horizon. Note, however, that not only is an NLP solution to the MIOP obtained on nearly every guidance cycle, but the state path constraints remain below their maximum limits for all $1000$ guidance trajectories because of the inclusion of the margin-inducing terms in the MIOP objective functional. 

\begin{table}[htbp]
\centering
\caption{Maximum heating rate and dynamic pressure results using the MIOG method (margin amounts in parentheses).}
\label{tab:max-path-cons-MIOG}

\centering

\begin{tabular}{lcccc}
\hline\hline
& \multicolumn{2}{c}{IPOPT}
& \multicolumn{2}{c}{SPRNLP} \\
\cline{2-3}\cline{4-5}
Statistic
& $\max(\dot{Q})$ & $\max(q)$
& $\max(\dot{Q})$ & $\max(q)$ \\
& (MW/m$^2$) & (kPa)
& (MW/m$^2$) & (kPa) \\
\hline
Reference
& $0.714\;(0.136)$ & $13.608\;(2.767)$
& $0.714\;(0.136)$ & $13.608\;(2.767)$ \\
Mean
& $0.713\;(0.137)$ & $12.894\;(3.481)$
& $0.717\;(0.133)$ & $13.608\;(2.767)$ \\
Max
& $0.772\;(0.0783)$ & $15.294\;(1.081)$
& $0.787\;(0.0632)$ & $15.843\;(0.532)$ \\
Std. Dev.
& $0.0207\;(0.0207)$ & $0.699\;(0.699)$
& $0.0232\;(0.0232)$ & $0.790\;(0.790)$ \\
\hline\hline
\end{tabular}

\end{table}

\begin{table}[htbp]
\centering
\caption{MIOG guidance cycles where an NLP solution was not obtained.}
\label{tab:feas-results}

\begin{tabular}{lcc}
\hline\hline
Statistic & IPOPT & SPRNLP \\
\hline
Mean      & $0.89$ & $0.96$ \\
Max       & $3$ & $3$ \\
Std. Dev. & $0.63$ & $0.58$ \\
\hline\hline
\end{tabular}
\end{table}

\subsubsection{Behavior of MIOG Guidance Trajectories}\label{sect:rlv-behavior-miog}

The second key feature of the MIOG method is that, not only is margin gained in the guidance trajectories, but the guidance trajectories also exhibit desired characteristics with only a small loss in performance relative to the problem where path constraints are enforced. Specifically, the objective functional of the path-constrained optimal control problem was designed to produce solutions with smooth control, reduce phugoid oscillations, and increase the terminal crossrange.  Each of these desired characteristics exhibited by the MIOG method is discussed in the remainder of this section. 

\begin{figure}[htb]
    \centering

    \begin{minipage}[t]{0.48\textwidth}
        \begin{subfigure}[t]{\linewidth}
            \centering
            \includegraphics[width=\linewidth]{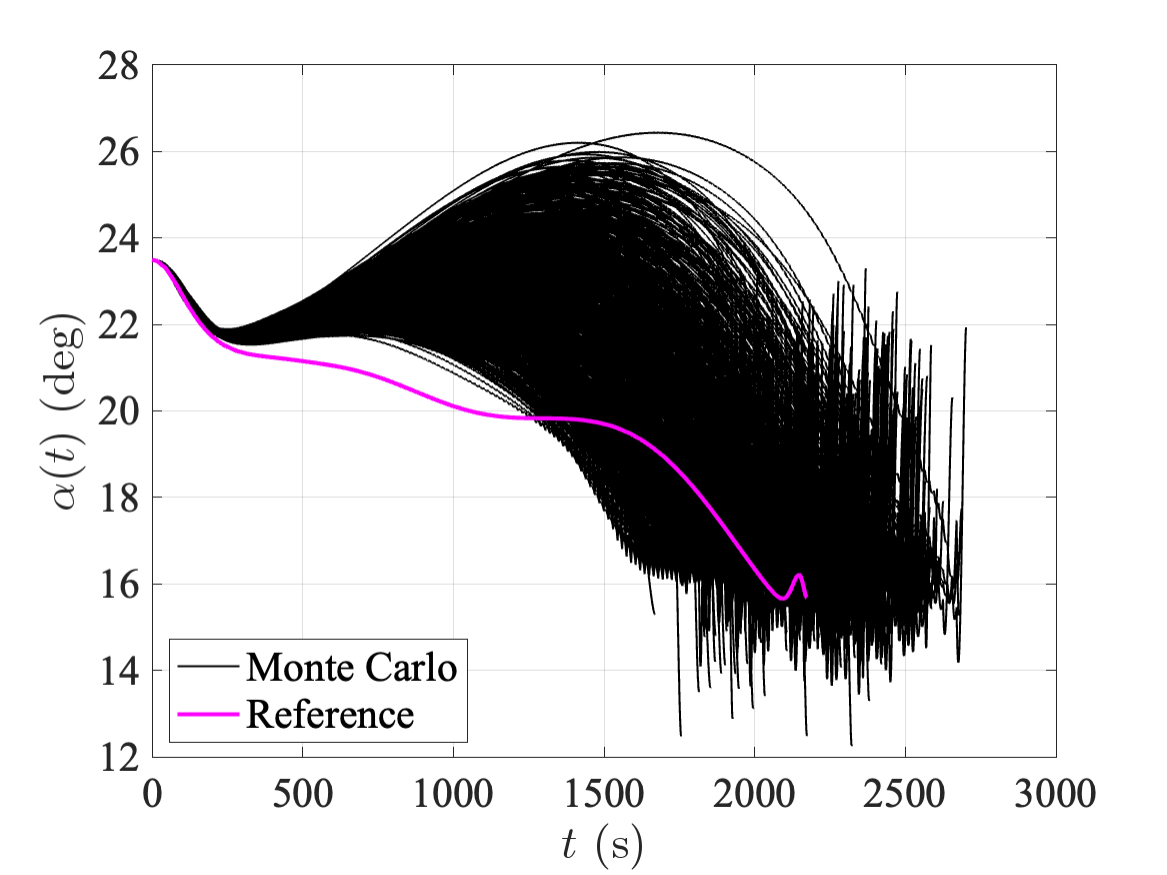}
            \caption{$\alpha(t)$ vs.~$t$ using IPOPT.}
            \label{fig:alpha-IPOPT}
        \end{subfigure}

        \vspace{1em}

        \begin{subfigure}[t]{\linewidth}
            \centering
            \includegraphics[width=\linewidth]{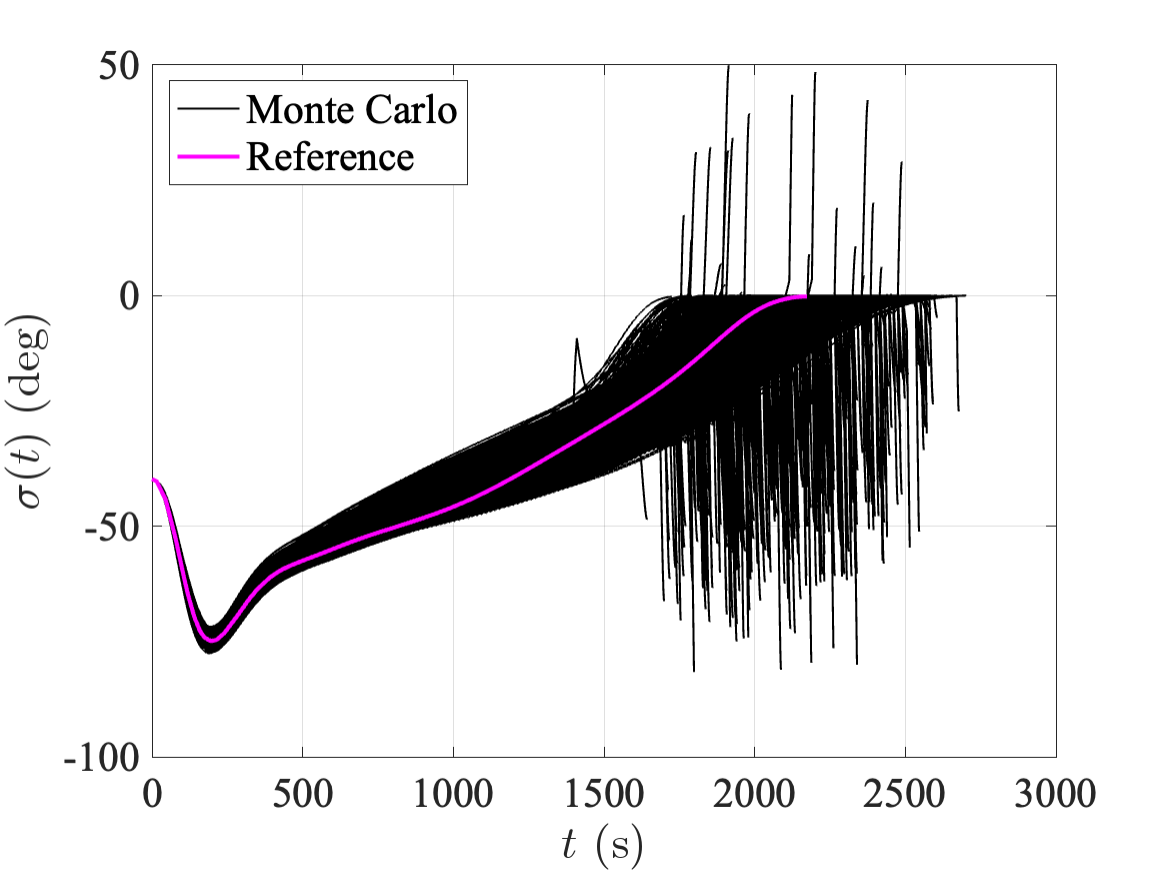}
            \caption{ $\sigma(t)$ vs.~$t$ using IPOPT.}
            \label{fig:bank-IPOPT}
        \end{subfigure}
    \end{minipage}
    \hfill
    \begin{minipage}[t]{0.48\textwidth}
        \begin{subfigure}[t]{\linewidth}
            \centering
            \includegraphics[width=\linewidth]{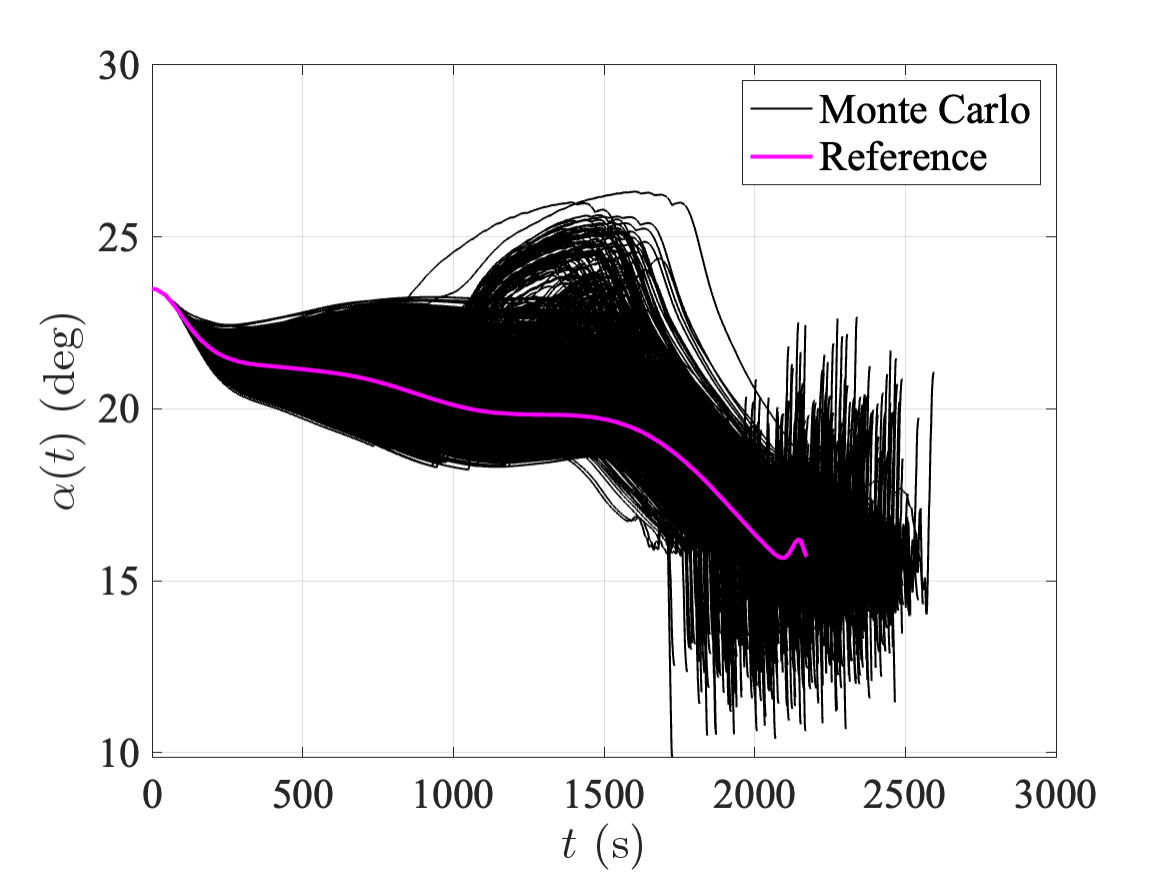}
            \caption{ $\alpha(t)$ vs.~$t$ using SPRNLP.}
            \label{fig:alpha-SOS}
        \end{subfigure}

        \vspace{1em}

        \begin{subfigure}[t]{\linewidth}
            \centering
            \includegraphics[width=\linewidth]{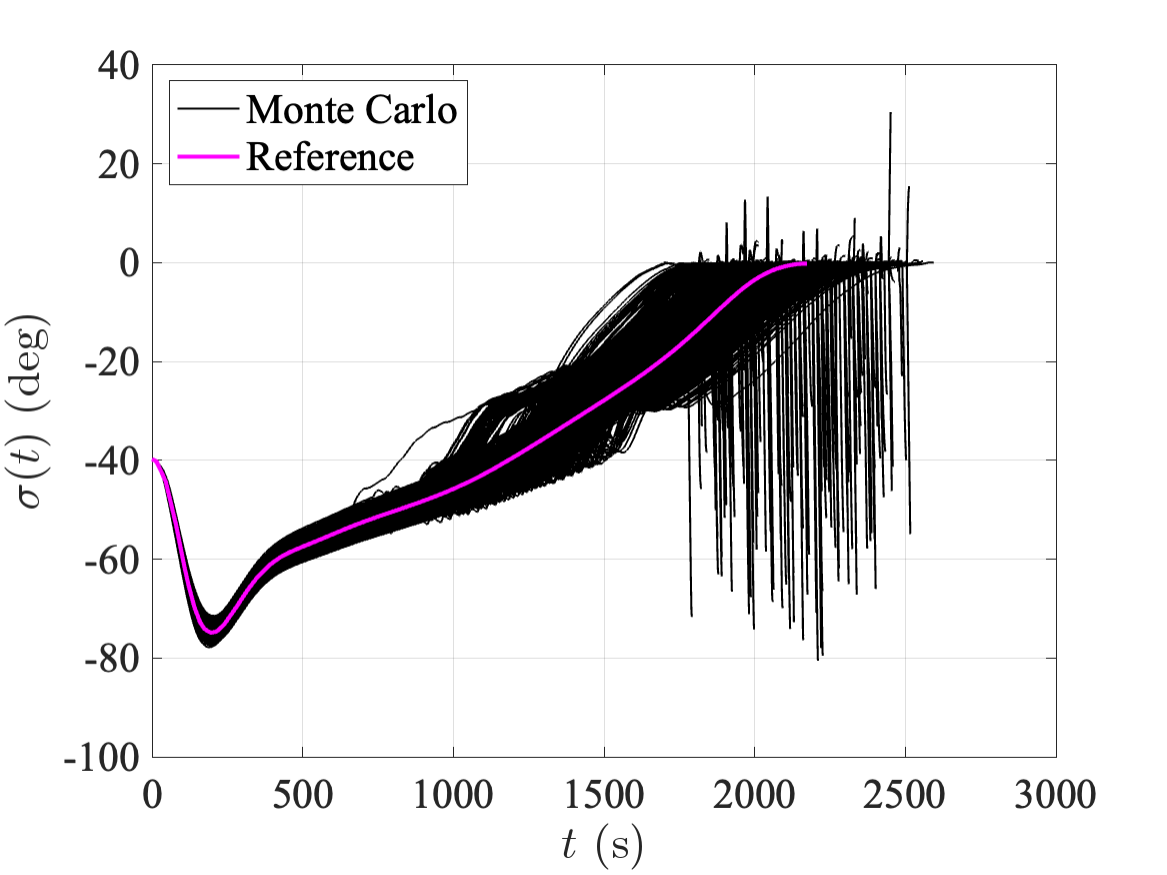}
            \caption{ $\sigma(t) $ vs.~$t$ using SPRNLP.}
            \label{fig:bank-SOS}
        \end{subfigure}
    \end{minipage}
    \caption{Angle of attack and bank angle obtained using the MIOG method with IPOPT and SPRNLP.}
    \label{fig:alpha-sigma-MIOG}
\end{figure}

\begin{figure}[htb]
    \centering    \includegraphics[width=0.5\linewidth]{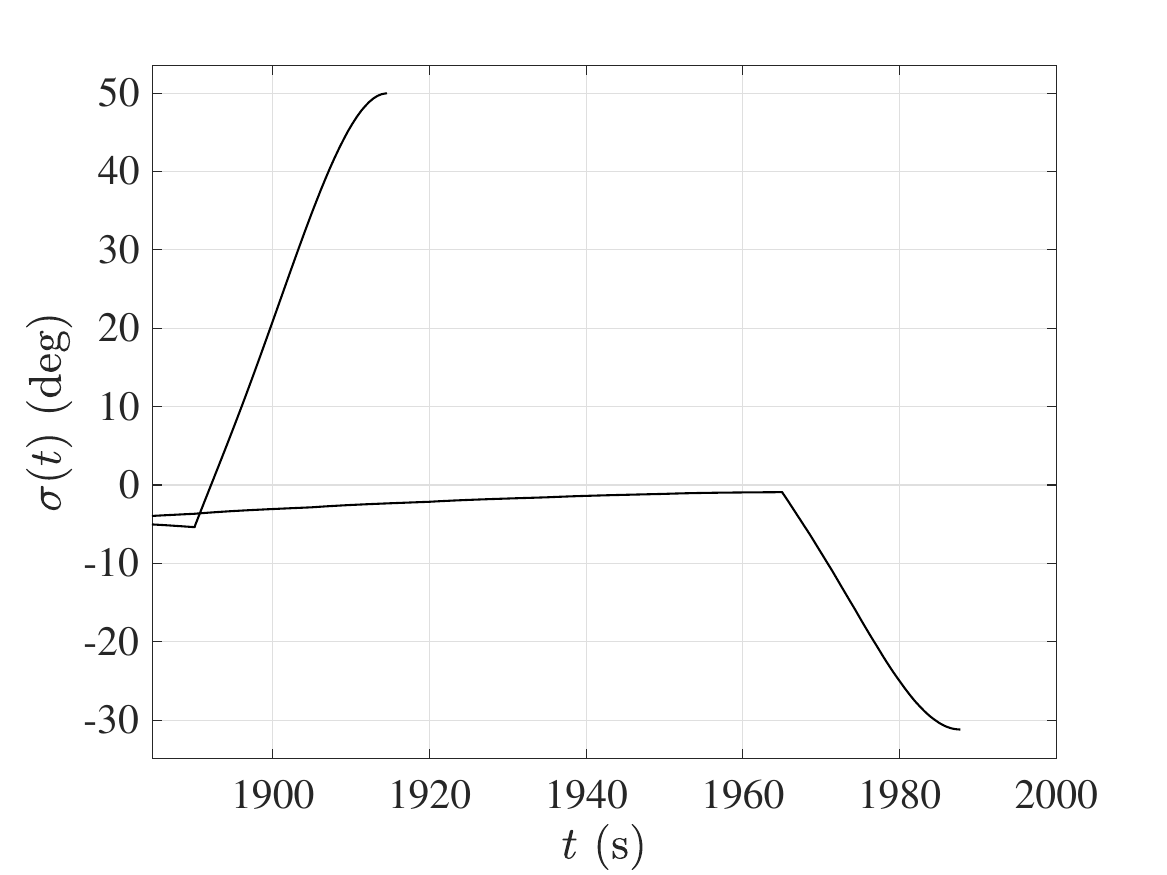}
    \caption{Enlarged view of two sample MIOG (with IPOPT) bank angle trajectories near terminal time.}
    \label{fig:bank-zoomed}
\end{figure}

With regard to producing a smooth control, Figs.~\ref{fig:alpha-IPOPT}--\ref{fig:bank-SOS} show that the MIOP reference angle of attack and bank angle are smooth.  Next, the angle of attack and bank angle on the guidance trajectories obtained using the MIOG method are also shown in Figs.~\ref{fig:alpha-IPOPT}--\ref{fig:bank-SOS}. It is noted that while in Figs.~\ref{fig:alpha-IPOPT}--\ref{fig:bank-SOS} $(\alpha,\sigma)$ appear to change rapidly near the end of the trajectories, this drastic appearance is a result of the axes scaling. For clarity, Fig.~\ref{fig:bank-zoomed} provides an enlarged view of two sampled bank angles near the final time, for which it can be seen that these maneuvers satisfy the control rate limits. It is noted, however, that the control changes more quickly towards the end of the guidance trajectories due to the need to satisfy the terminal constraints.

\begin{figure}[htb]
    \centering

    \begin{minipage}[t]{0.48\textwidth}
        \begin{subfigure}[t]{\linewidth}
            \centering
            \includegraphics[width=\linewidth]{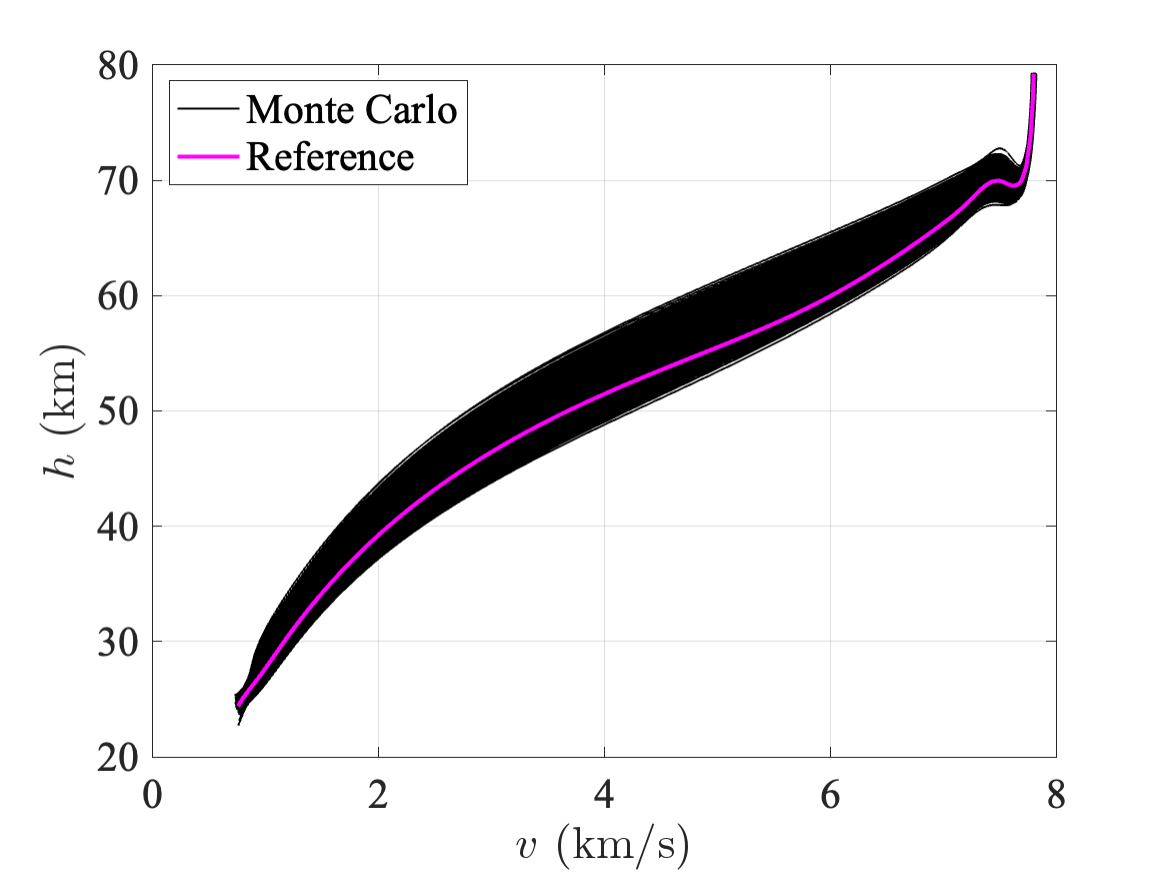}
            \caption{ $h$ vs.~$v$ using IPOPT.}
            \label{fig:alt-v-speed-IPOPT}
        \end{subfigure}

        \vspace{1em}

        \begin{subfigure}[t]{\linewidth}
            \centering
            \includegraphics[width=\linewidth]{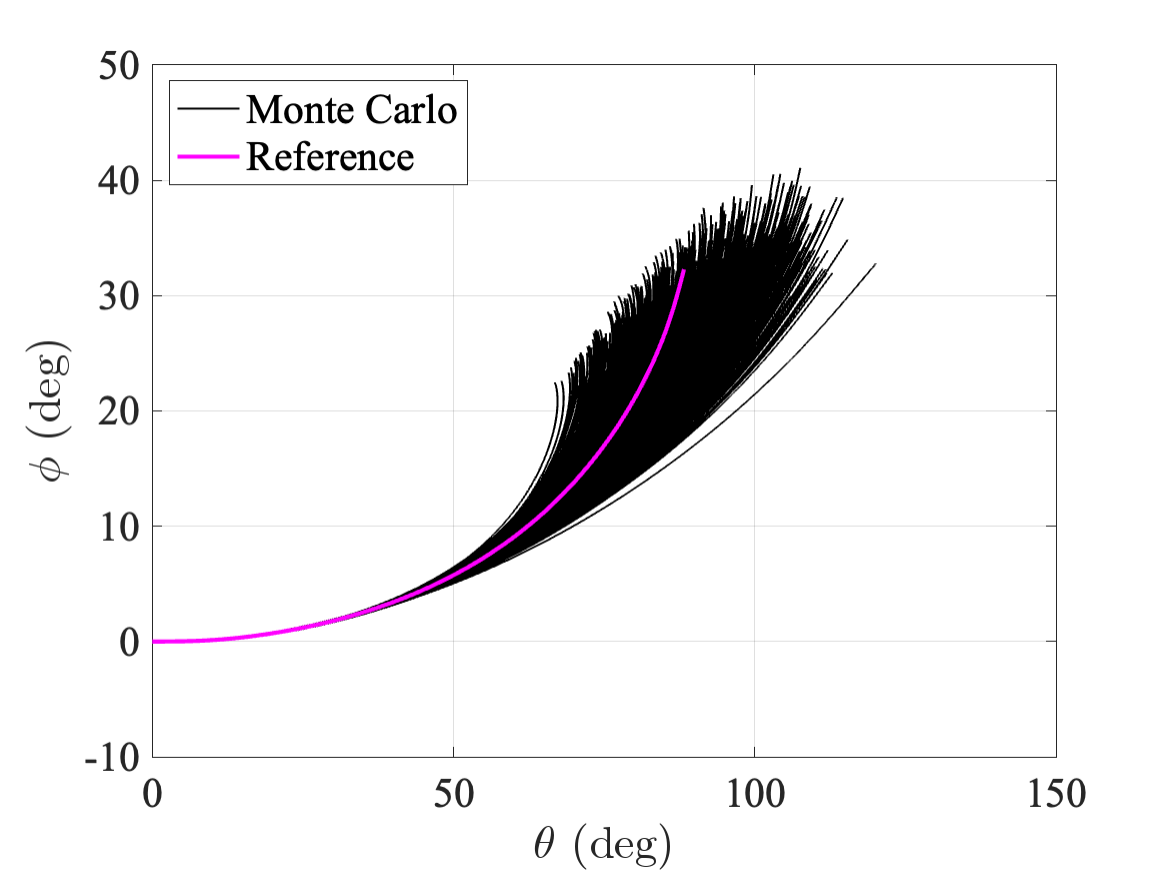}
            \caption{ $\phi$ vs.~$\theta$ using IPOPT.}
            \label{fig:lon-lat-IPOPT}
        \end{subfigure}

         \vspace{1em}

        \begin{subfigure}[t]{\linewidth}
            \centering
            \includegraphics[width=\linewidth]{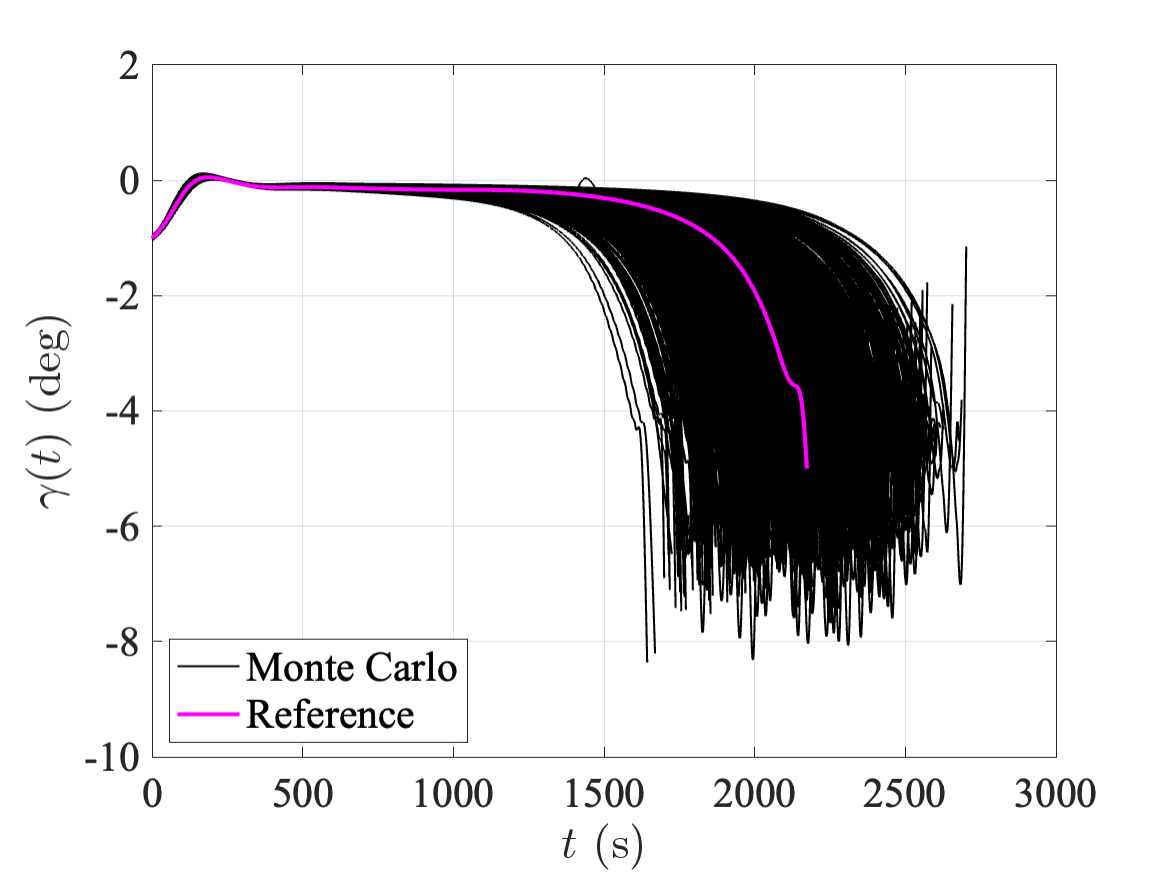}
            \caption{ $\gamma(t)$ vs.~$t$ using IPOPT.}
            \label{fig:fpa-IPOPT}
        \end{subfigure}
    \end{minipage}
    \hfill
    \begin{minipage}[t]{0.48\textwidth}
        \begin{subfigure}[t]{\linewidth}
            \centering
            \includegraphics[width=\linewidth]{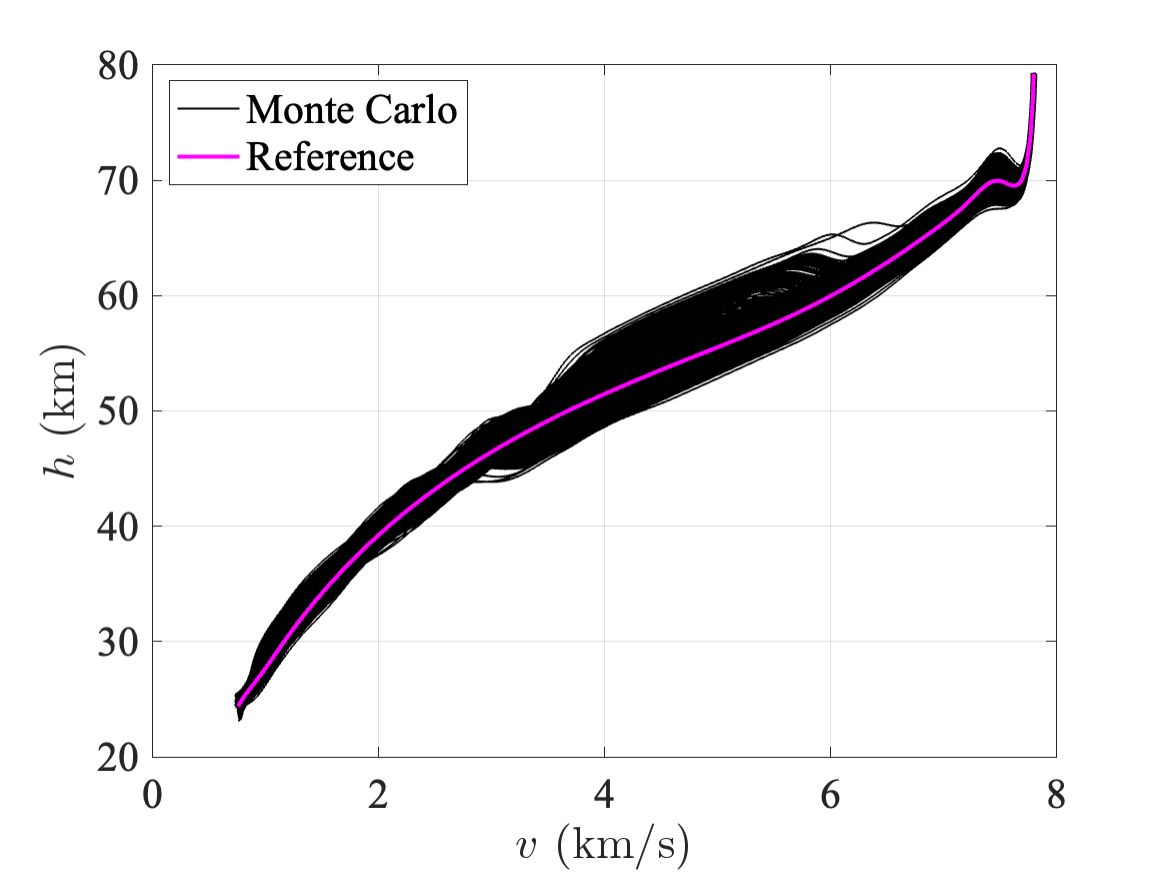}
            \caption{ $h$ vs.~$v$ using SPRNLP.}
            \label{fig:alt-v-SOS}
        \end{subfigure}

        \vspace{1em}

        \begin{subfigure}[t]{\linewidth}
            \centering
            \includegraphics[width=\linewidth]{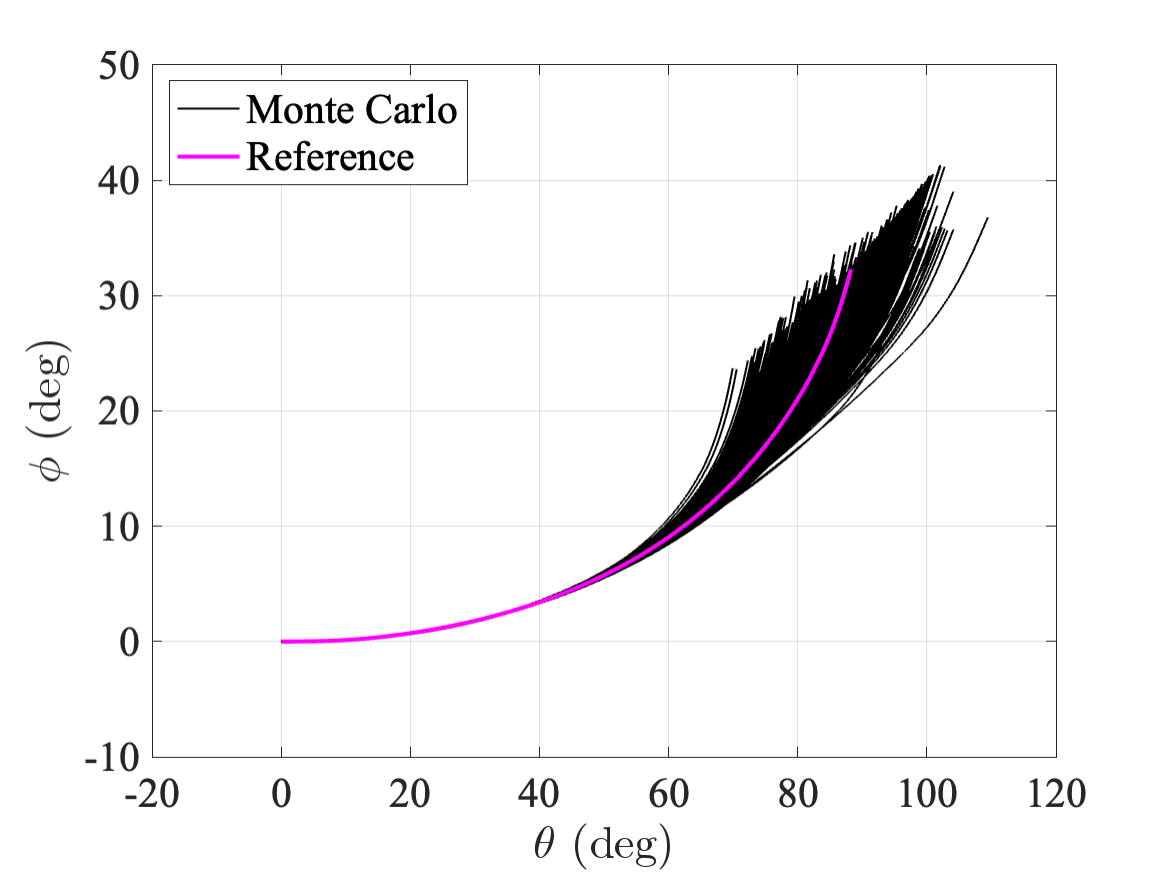}
            \caption{ $\phi $ vs.~$\theta$ using SPRNLP.}
            \label{fig:lon-lat-SOS}
        \end{subfigure}

         \vspace{1em}

        \begin{subfigure}[t]{\linewidth}
            \centering
            \includegraphics[width=\linewidth]{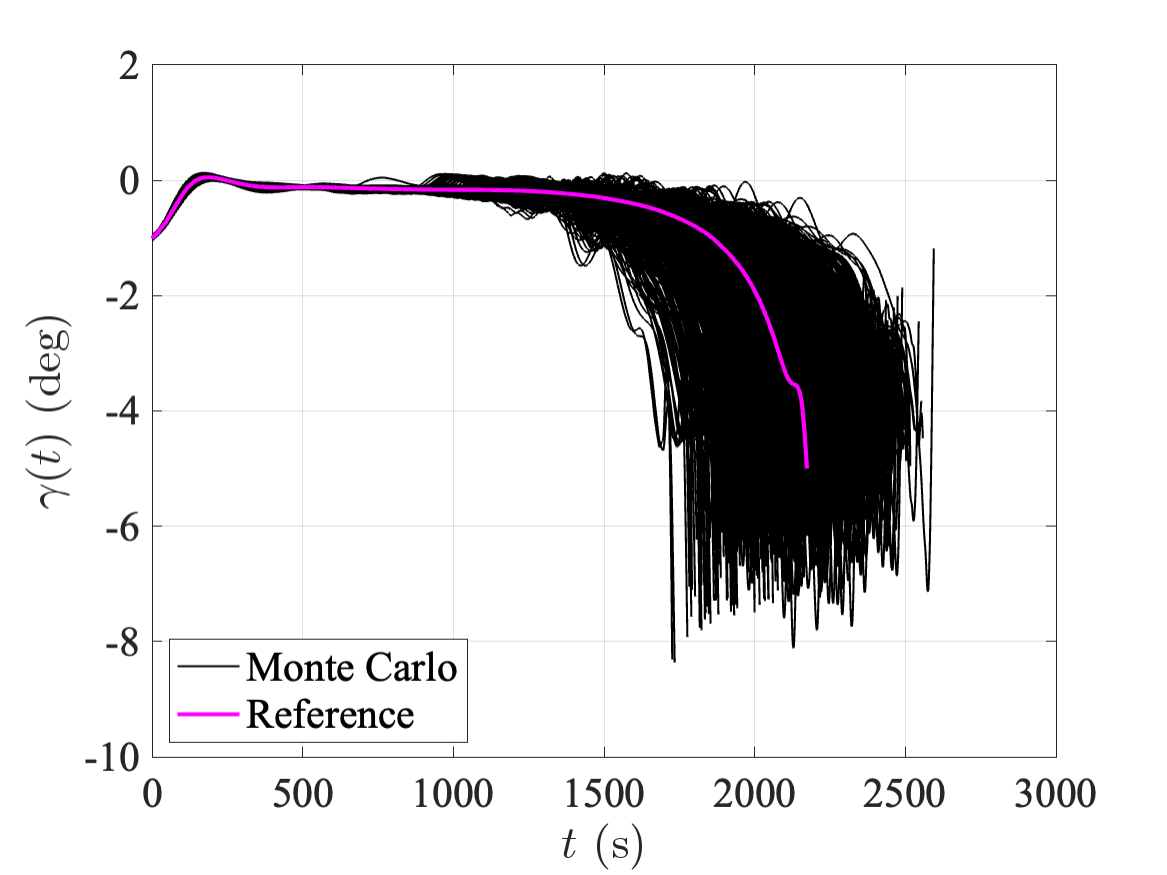}
            \caption{ $\gamma(t) $ vs.~$t$ using SPRNLP.}
            \label{fig:fpa-SOS}
        \end{subfigure}
    \end{minipage}

    \caption{State obtained using the MIOG method with IPOPT and SPRNLP.}
    \label{fig:states-MIOG}
\end{figure}

With regard to phugoid oscillations, Figs.~\ref{fig:alt-v-speed-IPOPT}--\ref{fig:fpa-SOS} show that the MIOP reference altitude and flight path angle do not contain phugoid oscillations. In addition, the MIOG guidance trajectories behave similarly to the MIOP reference. It is noted that, using SPRNLP, the flight path angle exhibits small fluctuations, but these fluctuations do not have a significant effect on the behavior of the guidance trajectory. 

With regard to increasing crossrange during entry, Table \ref{tab:terminal-lat} shows the statistics on the terminal latitude, $\phi(t_f)$, using the MIOG method. For comparison, the COP reference produced a terminal latitude of $33.08$~deg using IPOPT and $33.09$~deg using SPRNLP.  It is seen that the terminal latitude of the MIOP reference solution is $2.41\%$ lower using IPOPT and $2.44\%$ lower using SPRNLP relative to the COP reference solution.  Next, the mean terminal latitudes using the MIOG method are $8.13\%$ lower using IPOPT and $2.63 \%$ lower using SPRNLP relative to the COP reference. The relatively small reduction in terminal latitude using the MIOG method is offset by a significant increase in path constraint margin.

\begin{table}[!htbp]
\centering
\caption{Terminal latitude results, $\phi(t_f)$~(deg), using the MIOG method.}
\label{tab:terminal-lat}

\begin{tabular}{lcc}
\hline\hline
Statistic & IPOPT & SPRNLP \\
\hline
Reference & $32.28$ & $32.28$ \\
Mean      & $30.39$ & $32.22$ \\
Max       & $41.08$ & $41.32$ \\
Min       & $22.22$ & $23.38$ \\
Std. Dev. & $3.80$ & $3.76$ \\
\hline\hline
\end{tabular}

\end{table}
\subsubsection{Computational Performance Results}\label{sect:rlv-cpu-times-miog}

Another feature of the MIOG method is that solutions are obtained quickly.  In particular, pairing the MIOG method with LGR collocation and full Newton NLP solvers (in this case IPOPT and SPRNLP) leads to a sparse NLP that can be solved efficiently.  This computational efficiency becomes more pronounced as the guidance cycles evolve because the NLP becomes smaller due to the mesh truncation and remapping described in Ref.~\cite{Dennis2019} and reviewed in Section \ref{sect:mesh-remapping}. Moreover, the guidance updates are generated quickly because the initial guess from the truncated and remapped trajectory provides an excellent initial guess for the NLP solvers.  Table~\ref{tab:MIOG-computational} shows the computational times (CPU times) achieved by the MIOG method.  It is seen that the mean CPU time required to obtain a solution is less than one second using either NLP solver which is more than an order of magnitude less than the guidance cycle duration $\Delta T = 15$~s. In particular, SPRNLP has an exceptionally low mean CPU time of $0.0428$~s.  This extremely fast SPRNLP solve time is due to the fact that SPRNLP employs a full Newton {\em active-set} sequential quadratic programming (SQP) method which takes full advantage of a good initial guess. Because the truncated and remapped solution from the previous remaining  horizon is the initial guess for the next remaining horizon, SPRNLP converges in a small number of iterations.  While the IPOPT CPU times are also small, they are not quite as small as those obtained using SPRNLP because IPOPT employs an interior point method which does not take full advantage of a good initial guess. Consequently, despite having a good initial guess, IPOPT takes significantly more iterations to converge compared with SPRNLP.  It is noted, however, that IPOPT still converges relatively quickly because IPOPT is a computationally efficient solver.

\begin{table}[!htbp]
\centering
\caption{Computational results (CPU time in seconds) using the MIOG method.}
\label{tab:MIOG-computational}

\centering

\begin{tabular}{lcccc}
\hline\hline
& \multicolumn{2}{c}{IPOPT}
& \multicolumn{2}{c}{SPRNLP} \\
\cline{2-3}\cline{4-5}
Statistic
& Mean CPU & Max CPU
& Mean CPU & Max CPU \\
\hline
Mean
& $0.176$ & $0.664$
& $ 0.0428$ & $ 0.441$ \\
Max
& $0.726$ & $5.052$
& $0.222$ & $9.256$ \\
Std. Dev.
& $0.0520$ & $0.439$
& $ 0.0212$ & $0.887$ \\
\hline\hline
\end{tabular}


\end{table}

\subsubsection{Accuracy of MIOG Method}\label{rlv-accuracy-miog}

\begin{figure}[h]
    \centering

    \begin{minipage}[t]{0.48\textwidth}
        \begin{subfigure}[t]{\linewidth}
            \centering
            \includegraphics[width=\linewidth]{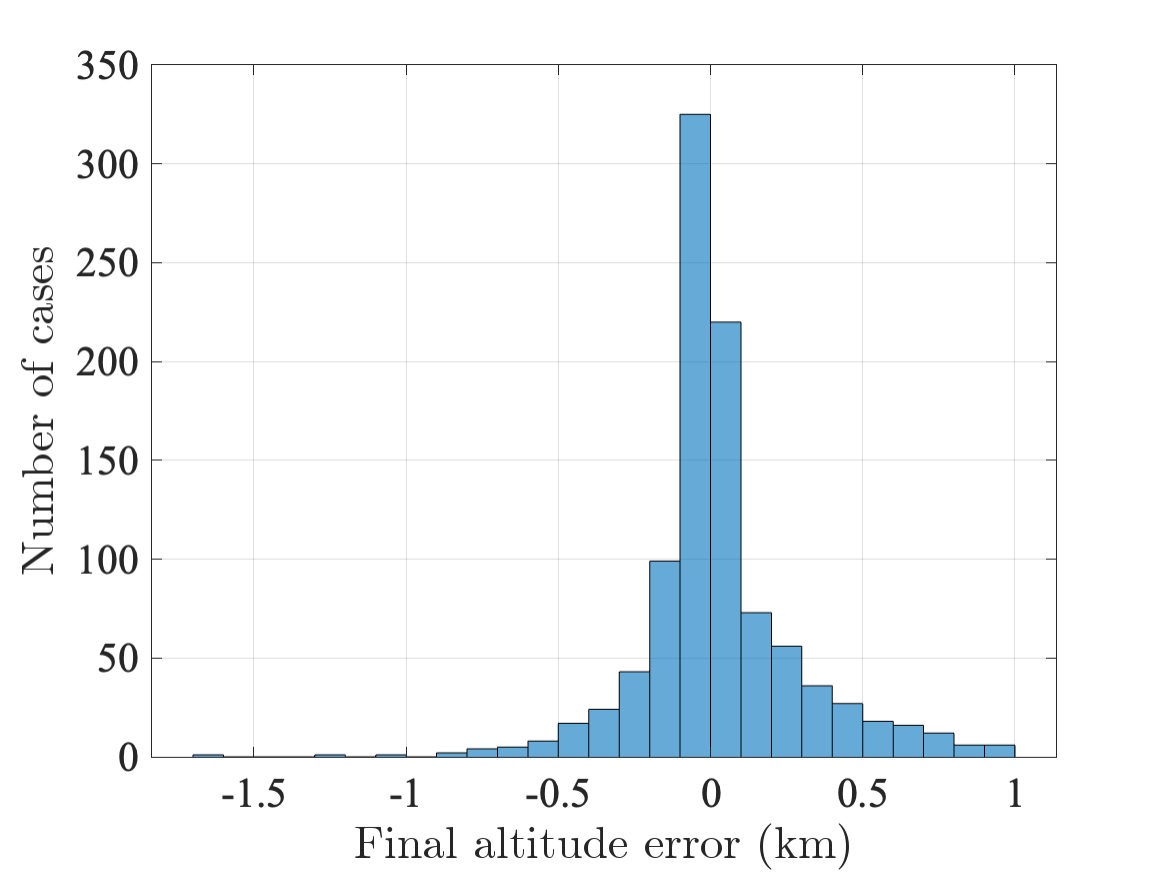}
            \caption{Terminal altitude error histogram using IPOPT.}
            \label{fig:alt-histogram-IPOPT}
        \end{subfigure}

        \vspace{1em}

        \begin{subfigure}[t]{\linewidth}
            \centering
            \includegraphics[width=\linewidth]{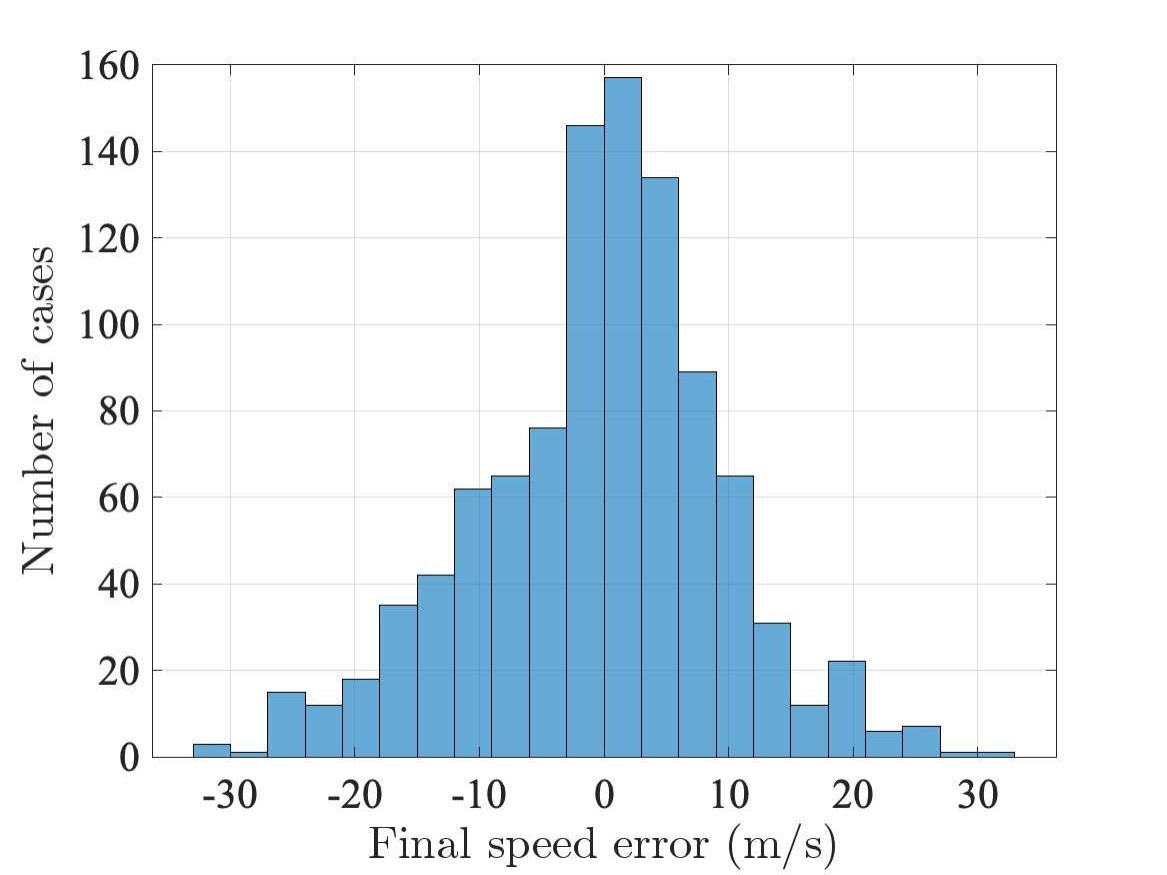}
            \caption{Terminal speed error histogram using IPOPT.}
            \label{fig:speed-histogram-IPOPT}
        \end{subfigure}

         \vspace{1em}

        \begin{subfigure}[t]{\linewidth}
            \centering
            \includegraphics[width=\linewidth]{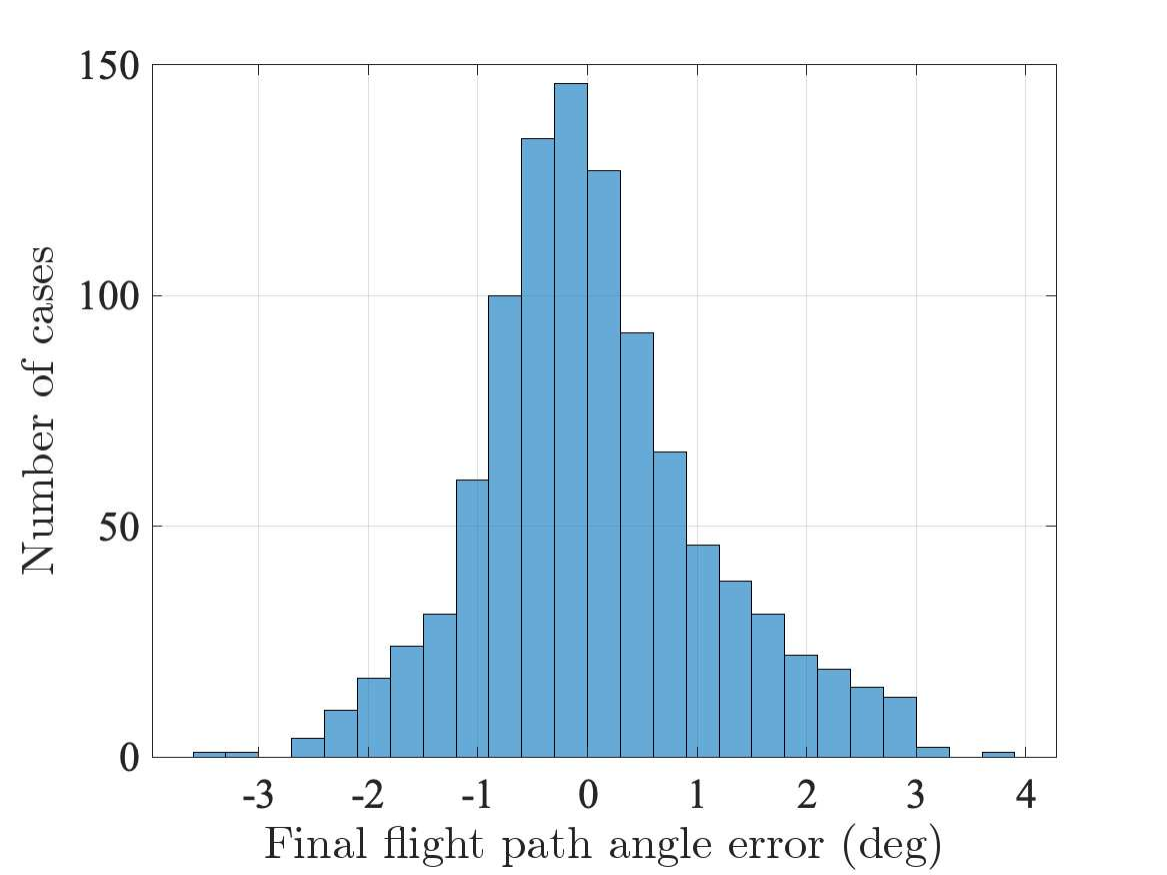}
            \caption{Terminal flight path angle error histogram using IPOPT.}
            \label{fig:fpa-histogram-IPOPT}
        \end{subfigure}
    \end{minipage}
    \hfill
    \begin{minipage}[t]{0.48\textwidth}
        \begin{subfigure}[t]{\linewidth}
            \centering
            \includegraphics[width=\linewidth]{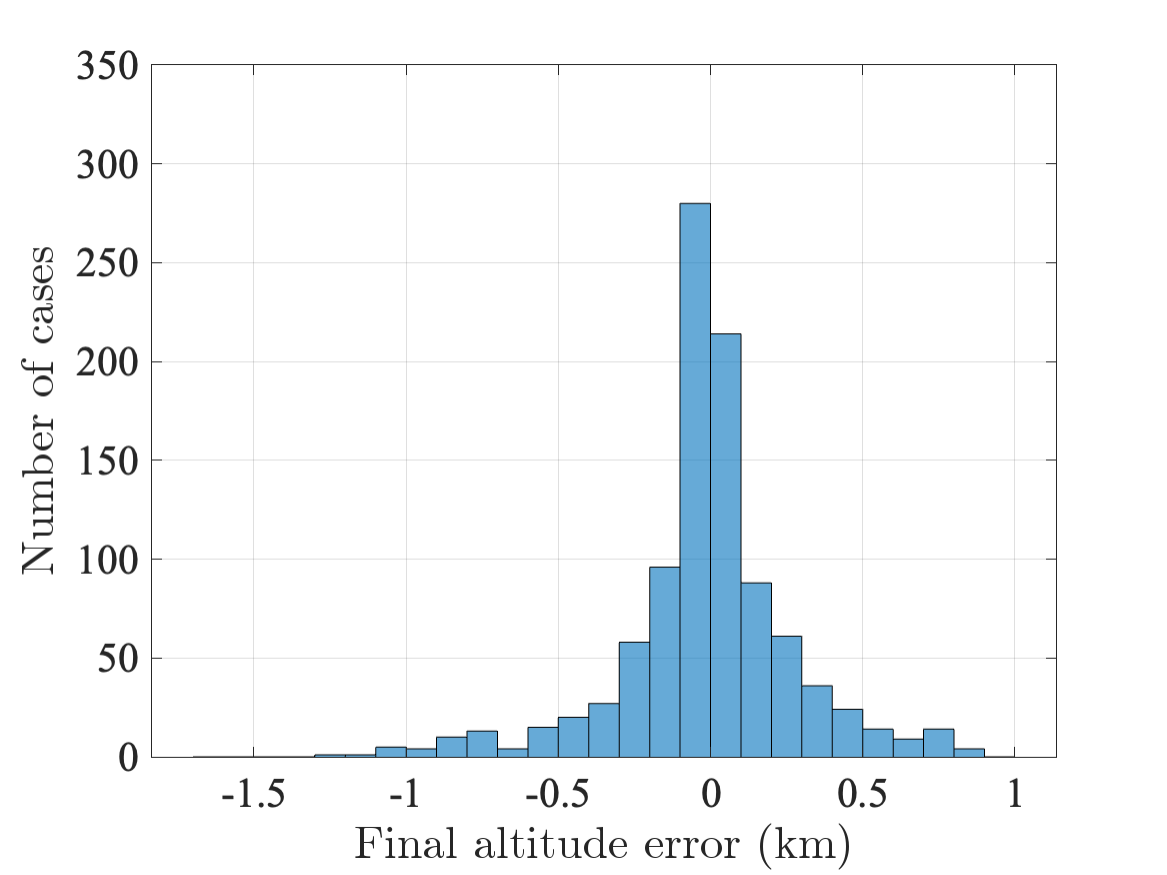}
            \caption{Terminal altitude error histogram using SPRNLP.}
            \label{fig:alt-histogram-SOS}
        \end{subfigure}

        \vspace{1em}

        \begin{subfigure}[t]{\linewidth}
            \centering
            \includegraphics[width=\linewidth]{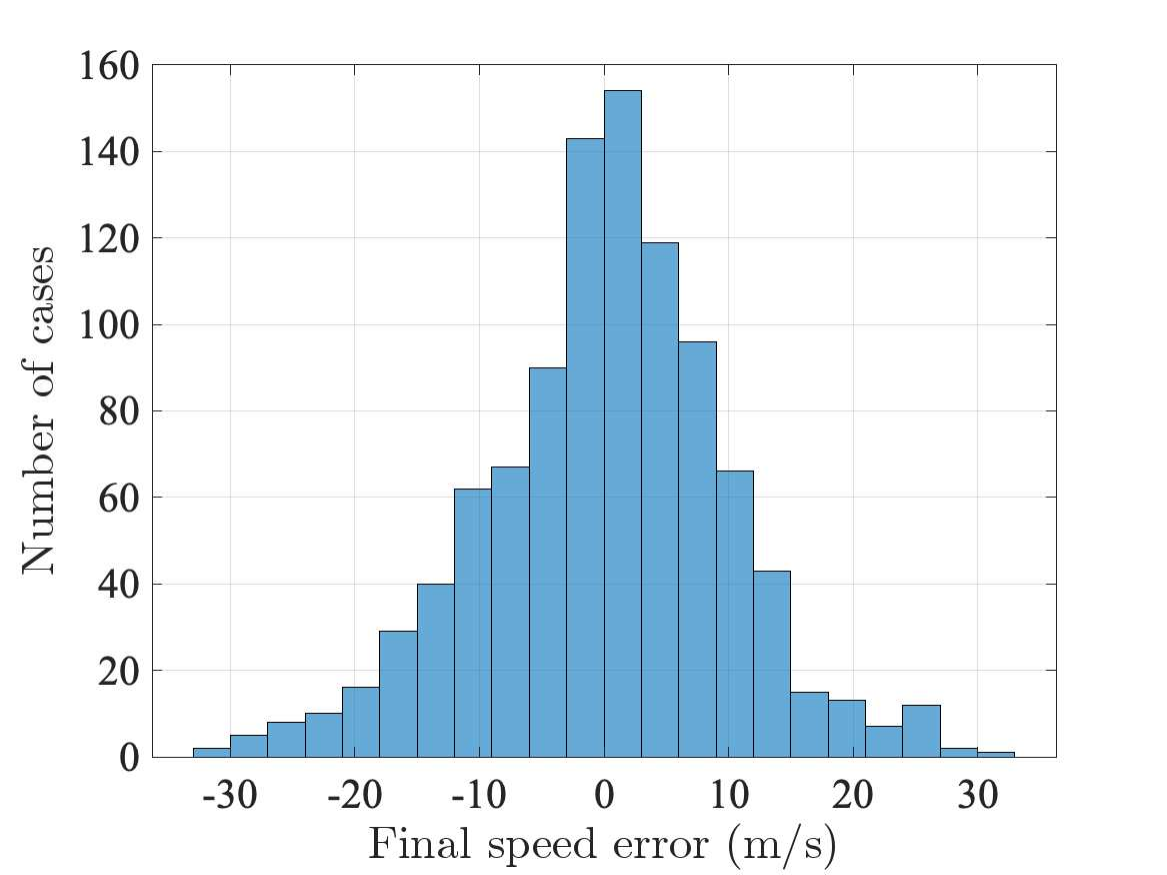}
            \caption{Terminal speed error histogram using SPRNLP.}
            \label{fig:speed-histrogram-SOS}
        \end{subfigure}

         \vspace{1em}

        \begin{subfigure}[t]{\linewidth}
            \centering
            \includegraphics[width=\linewidth]{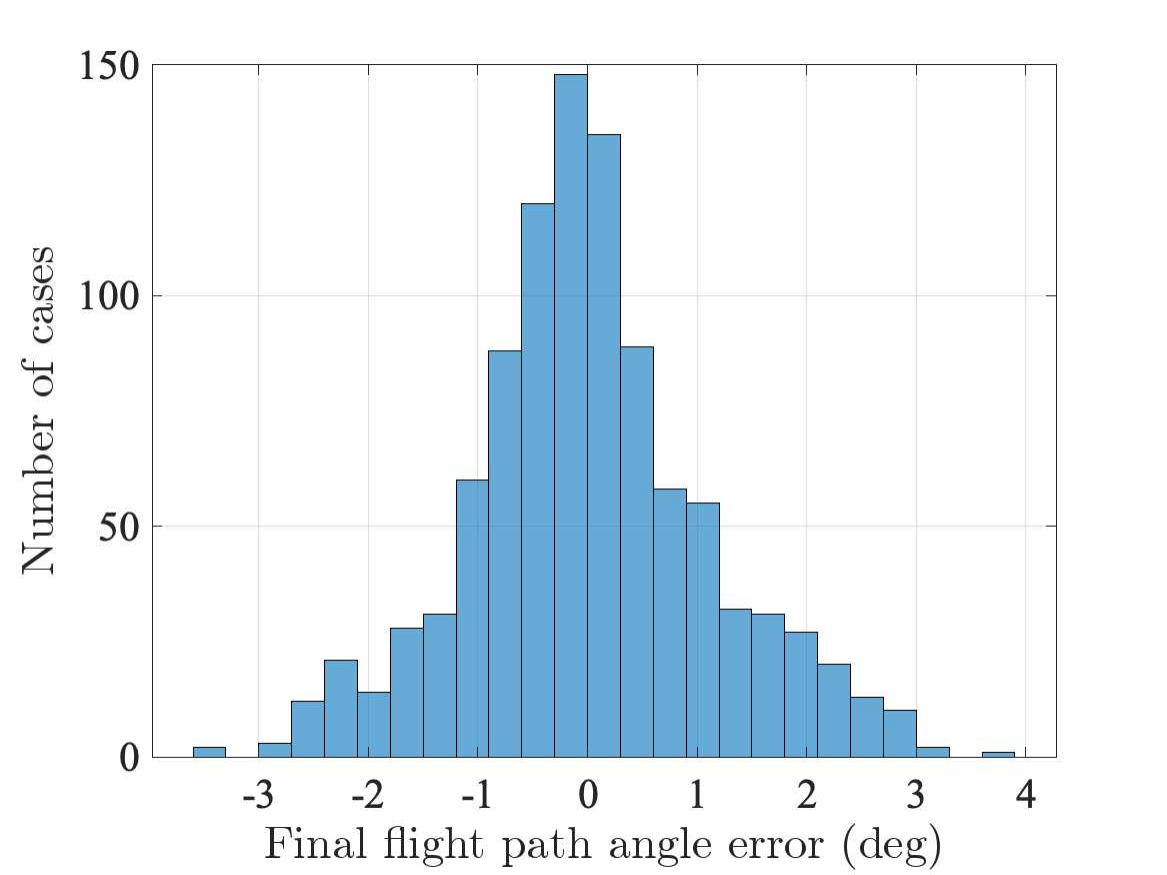}
            \caption{Terminal flight path angle error histogram using SPRNLP.}
            \label{fig:fpa-histogram-SOS}
        \end{subfigure}
    \end{minipage}

    \caption{Terminal state error histograms obtained using the MIOG method with IPOPT and SPRNLP.}
    \label{fig:hist-MIOG}
\end{figure}

Finally, the last key feature of the MIOG method is the accuracy with which the terminal constraints are attained. The mean and maximum values of the absolute terminal state errors are given in Table~\ref{tab:terminal-state-errors} alongside the standard deviations of the terminal state errors.  It is seen in Table~\ref{tab:terminal-state-errors} and the histograms in Fig.~\ref{fig:hist-MIOG} that the mean absolute error using the MIOG method is under $200$~m in altitude, under $10$~m/s in speed, and under $1$~deg in flight path angle using either NLP solver. Furthermore, the maximum absolute terminal errors of the altitude, speed, and flight path angle are less than $2000$~m, $40$~m/s, and $5$~deg, respectively. These results demonstrate that the MIOG method is able to compensate for uncertainties.

\begin{table}[!htbp]
\centering
\caption{Terminal state errors (absolute) and standard deviations (signed) using the MIOG method.}
\label{tab:terminal-state-errors}

\begin{tabular}{lcccccc}
\hline\hline
& \multicolumn{3}{c}{IPOPT}
& \multicolumn{3}{c}{SPRNLP} \\
\cline{2-4}\cline{5-7}
Statistic
& $\delta h(t_f)$ & $\delta v(t_f)$ & $\delta\gamma(t_f)$
& $\delta h(t_f)$ & $\delta v(t_f)$ & $\delta\gamma(t_f)$ \\
& (m) & (m/s) & (deg)
& (m) & (m/s) & (deg) \\
\hline
Mean (Absolute)
& $166.2$ & $7.47$ & $0.80$
& $184.9$ & $7.45$ & $0.830$ \\
Max (Absolute)
& $1677$ & $31.53$ & $3.85$
& $1292$ & $32.49$ & $3.82$ \\
Standard Deviation (Signed)
& $260.1$ & $9.85$ & $1.05$
& $286.5$ & $9.81$ & $1.10$ \\
\hline\hline
\end{tabular}

\end{table}

\clearpage

\section{Summary of Results}\label{sect:summary}

The results presented in Section~\ref{sect:rlv-results-miog} highlight the following key aspects of the MIOG method developed in this paper. First, the MIOP reference induces significant margin in the heating rate and dynamic pressure.  Furthermore, the corresponding guidance trajectories maintain significant margin with a wide range of perturbations from the reference solution. In fact, margin is induced by the MIOG method in every guidance trajectory (that is, for every Monte Carlo trial) using the MIOG method. Second, the guidance trajectories obtained using the MIOG method retain a large portion of the performance relative to the solution of the original COP. In addition, the control obtained using the MIOG method produces guidance trajectories with desirable characteristics including mitigating phugoid oscillations and large terminal crossrange. The third key feature is that the MIOG method is computationally efficient, rapidly obtaining solutions.  In particular, the average computational time is under a half second using either IPOPT or SPRNLP. The fourth key feature is that an NLP solution to the MIOP is obtained on nearly every guidance cycle, ensuring that the control can be updated on nearly every guidance cycle.  Finally, the MIOG method produces guidance trajectories where the errors in the terminal state are small.

\section{Conclusions}\label{sect:conclusions}

A method has been developed for margin-inducing optimal guidance and control with application to atmospheric entry.  The key features of the method are the ability of the guidance trajectory to increase margin in state inequality path constraints while retaining a large portion of the performance of the original path-constrained optimal control problem. A key part of the method is to develop an objective functional that employs a smooth approximation of the maximum for a function of a continuous variable.  This smooth approximation is derived from a modification of the Kreisselmeier-Steinhauser function which is defined for a function over sampled data.  This smooth approximation is then used to formulate a margin-inducing optimal control problem where the objective function contains margin-inducing terms based on the smooth approximation of the maximum function.  This margin-inducing optimal control problem is then re-solved at specified guidance update times to produce a guidance trajectory that maintains good performance with respect to the original objective functional while maintaining margin in the state inequality path constraints.  The method is applied to a state inequality path-constrained reusable launch vehicle atmospheric entry problem where it is shown to be computationally efficient and significantly increases margin while maintaining performance and accuracy.  The numerical results demonstrate that the method maintains high performance while avoiding path constraint violations in the guidance trajectories.

\section*{Acknowledgments}

The authors acknowledge support for this research from the U.S.~Office of Naval Research under grant N00014-22-1-2397, from the U.S.~National Science Foundation under grant CMMI-2031213, and from the U.S.~Air Force Research Laboratory under grant FA8651-24-1-0004.

\renewcommand{\baselinestretch}{1}
\normalsize\normalfont

\bibliographystyle{aiaa}

\end{document}